\documentclass{article}

\usepackage{PRIMEarxiv}

\usepackage[utf8]{inputenc} 
\usepackage[T1]{fontenc}    
\usepackage{hyperref}       
\usepackage{url}            
\usepackage[table]{xcolor}
\usepackage{booktabs}       
\usepackage{amsfonts}       
\usepackage{nicefrac}       
\usepackage{microtype}      
\usepackage{lipsum}
\usepackage{fancyhdr}       
\usepackage{graphicx}       
\graphicspath{{media/}}     

\usepackage{hyperref}
\usepackage{url}            
\usepackage{booktabs}       
\usepackage{amsfonts}       
\usepackage{nicefrac}       
\usepackage{microtype}      
\usepackage{lipsum}		
\usepackage{natbib}
\usepackage{doi}
\usepackage{amsmath, amssymb}
\usepackage{soul, xcolor}
\usepackage{caption}
\usepackage{subcaption}
\usepackage{bm}
\usepackage{amsmath,amsfonts,amscd,amssymb,bm}
\usepackage{tikz}

\usepackage{pdflscape} 

\usepackage{soul, xcolor}
\usepackage{mathabx}
\usepackage{algorithm}
\usepackage{algpseudocode}
\usepackage{mathtools}

\newcommand{\R}{\mathbb{R}}

\newcommand{\enc}{\boldsymbol{\phi}}
\newcommand{\dec}{\boldsymbol{\psi}}
\newcommand{\xirom}{\boldsymbol{\xi}}

\title{Spatio-temporal Latent Denoising Diffusion Probabilistic Models for Reduced-order Modeling of Parametrized Dynamical Systems}

\author{
  Michiel Nikken \\
  Department of Applied Mathematics \\
  University of Twente \\
  Enschede, The Netherlands \\
  \texttt{m.s.nikken@utwente.nl} \\
  \and
  Nicolò Botteghi \\
  Modeling and Scientific Computing Laboratory \\
  Politecnico di Milano \\
  Milano, Italy \\
  \texttt{nicolo.botteghi@polimi.it} \\
  \and 
  Federico Califano \\
  Robotics and Mechatronics \\
  University of Twente \\
  Enschede, The Netherlands \\
  \texttt{f.califano@utwente.nl}
  \and
  Silke Glas \\
  Department of Applied Mathematics \\
  University of Twente \\
  Enschede, The Netherlands \\
  \texttt{s.m.glas@utwente.nl}
}

\begin{document}
\maketitle

\begin{abstract}
Many scientific problems require accurate modeling of complex physical phenomena, such as fluid dynamics or climate modeling. These phenomena often result in high-dimensional and thus computationally expensive computational models that can limit their application to real-time and multi-query problems. 
Model-order reduction (MOR) is an approach that seeks to approximate full-order models (FOMs) using reduced-order models (ROMs), trading a minor reduction in accuracy for a major reduction in computational cost. In this work, we propose a non-intrusive MOR method using generative machine learning by means of \emph{denoising diffusion probabilistic models} (DDPMs) for generating solutions of the dynamical systems under different instances of their parameters. Unlike conventional DDPMs, which often operate purely in the spatial domain, we aim to generate \emph{spatio-temporal} solutions to improve quality and temporal coherence. In addition, we embed the DDPM in a \emph{latent space} obtained by sequentially applying proper orthogonal decomposition and an autoencoder to reduce the data dimensionality and the computational cost of the DDPM. We test our approach on a parametrized 2D fluid flow around an obstacle. The numerical experiments demonstrate that our \emph{latent DDPM} can \emph{(i)} produce accurate and temporally coherent solutions, \emph{(ii)} achieve strong generalization capabilities to scenarios involving unseen parameter values, and \emph{(iii)} extrapolate in time beyond the training horizon.  
\end{abstract}

\section{Introduction}
%


Many scientific problems require accurate modeling of complex physical phenomena, which are often governed by nonlinear, time-dependent, and parameterized partial differential equations (PDEs) \cite{brunton2019data}. Examples arise in aerodynamics \cite{Paranjape2013}, energy systems \cite{ZHANG2022120081}, climate modeling \cite{Mengaldo2019}, chemical processes \cite{MARCICKI2013310}, and biomedical engineering \cite{Hawkins2012}, where predictive simulations are essential for design, monitoring, and decision making. For these types of problems, closed-form solutions are typically not available, so numerical approximations are used instead. The numerical solution of PDEs requires spatial and temporal discretization, yielding high-dimensional dynamical systems whose computational cost can be prohibitive \cite{quarteroni2015reduced}. These high-fidelity discretizations, commonly referred to as full-order models (FOMs), often involve thousands of degrees of freedom and demand large-scale computational resources. The computational cost of FOMs limits their application to real-time and multi-query problems 
such as optimal control \cite{ravindran2000reduced, troltzsch2024optimal}, uncertainty quantification \cite{sudret2000stochastic, galbally2010non}, and digital twins \cite{chakraborty2021role, willcox2024role}. Consequently, developing accurate yet computationally efficient surrogate models is crucial for such applications. Model-order reduction (MOR) techniques address this challenge by constructing reduced-order models (ROMs) that retain the essential dynamics of the original system while significantly reducing computational complexity.

Traditional projection-based MOR methods, like the POD-Galerkin method \cite{Aubry_Holmes_Lumley_Stone_1988}, obtain a ROM by projecting the FOM dynamics onto a low-dimensional linear space. However, in many scenarios, full knowledge of the FOM or even the underlying PDE is not available, e.g. when the FOM is part of proprietary simulation software. In that case, ROMs can still be constructed using data from the system. Such data-driven approaches are known as non-intrusive MOR. Notable non-intrusive MOR methods include operator inference \cite{peherstorfer2016data,guo2022bayesian}, sparse identification of reduced dynamics \cite{brunton2016discovering}, and dynamic mode decomposition \cite{schmid2010dynamic,GILL2026114718}.

In the last years, machine learning methods, especially using neural networks, have become increasingly popular for constructing non-intrusive ROMs. Firstly, they can be used to find efficient non-linear low-dimensional trial manifolds to represent the FOM state. This can mitigate limitations of linear approximation spaces, for instance, those given by the Kolmogorov $N$-width. The Kolmogorov $N$-width is the minimum worst-case approximation error of an $N$-dimensional linear approximation space. Particularly for transport-dominated problems, the Kolmogorov $N$-width may decay slowly with increasing $N$ \cite{Ohlberger2016Reduced}. This indicates that a large linear space is needed to find good approximations of the FOM. Since neural networks are general function approximators that do not require a priori assumptions on the (non-)linearity, they are effective tools to break this barrier. Secondly, machine learning methods are also used to model the dynamics in the reduced space, since in the non-intrusive case, these cannot be inferred from the FOM dynamics and even the structure may be unknown.

Examples of these machine learning methods include extensions to operator inference \cite{qian2020lift,SHARMA2024116865}, dynamic mode decomposition \cite{li2017extended,otto2019linearly}, and sparse identification of reduced dynamics \cite{champion2019data, bakarji2022discovering}. Other examples of non-intrusive MOR methods using machine learning are neural operators \cite{li2020fourier, kovachki2023neural, cao2023lno}, manifold learning with deep autoencoders \cite{fresca2021comprehensive,KimCWZ22}, Gaussian process and deep kernel for learning probabilistic dynamics \cite{guo2019data, botteghi2022deep, botteghi2024recurrentdeepkernellearning}, and generative models, such as generative adversarial networks, for learning generative dynamics \cite{kadeethum2021framework, kemna2023reduced, coscia2024generative}. 

Among the recently developed generative artificial intelligence (AI) approaches, denoising diffusion probabilistic models (DDPMs) \cite{Sohl-Dickstein2015DeepThermodynamics, Ho2020DenoisingModels} have emerged as a powerful framework to learn complex high-dimensional distributions and to generate high-quality data, such as images, audio, and molecular structures \cite{yang2023diffusion}. These DDPMs are based on a two-step probabilistic process: (i) a forward diffusion process that gradually adds noise to data, and (ii) a reverse process that reconstructs the original data by progressively removing the noise. Neural networks are used to learn the steps of the reverse process, such that repeated application of the denoising step recovers the original (noiseless) data from noisy samples. 

In the context of PDEs and scientific computing, DDPMs provide a flexible probabilistic framework for learning distributions over solution manifolds. By treating FOM solutions as samples from an unknown distribution, DDPMs can capture complex spatial correlations and nonlinear dynamics. Recent works have leveraged DDPMs for scientific computing problems, including autoregressive forecasting of PDE solutions \cite{price2023gencast, shysheya2024conditional, zhang2024xddpm}, data assimilation and reconstruction from partial measurements \cite{huang2024diffusionpde}, in combination with physics-informed objectives that enforce governing laws \cite{bastek2024physics, shan2024pird, shu2023physics}, and even generative offline control of PDEs \cite{wei2024generative}. Beyond PDE solvers, DDPMs have been successfully applied to molecular design \cite{xu2022geodiff}, microstructure reconstruction \cite{dureth2023conditional}, and metamaterial design \cite{bastek2023inverse}. Despite their success, directly applying DDPM to high-dimensional data still remains computationally demanding. This is because evaluating DDPMs require numerous forward passes through a neural network, of which the cost grows with the dimensionality of the data.


\begin{figure}[h!]
    \centering
    \includegraphics[width=\textwidth]{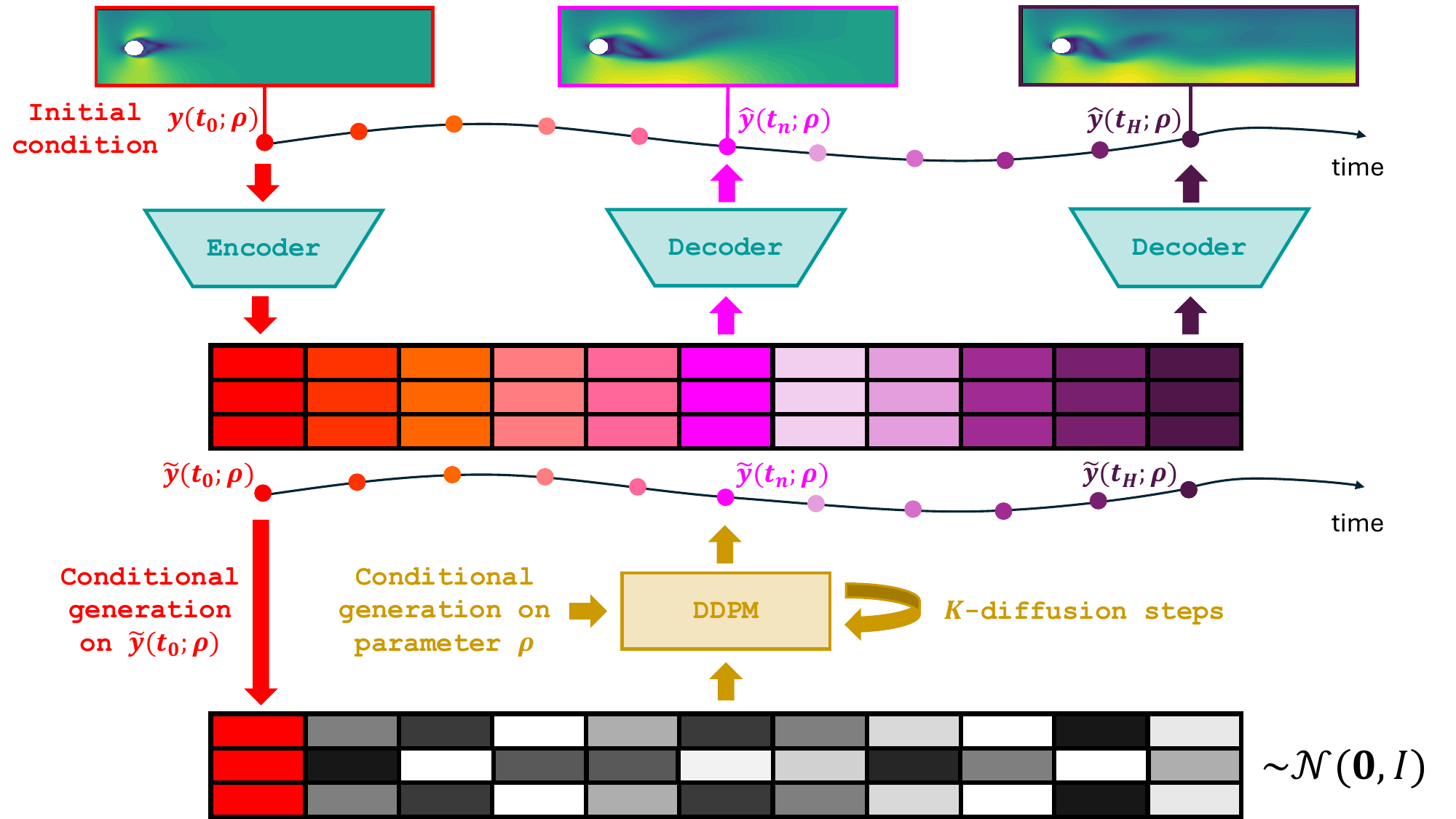}
    \caption{Latent denoising diffusion probabilistic model for generation of spatio-temporal solutions of parameterized dynamical systems. In the top-left, an initial condition is mapped to a reduced space in an encoding step that uses POD and an autoencoder. Next, $H$ noisy samples are sampled from a standard normal distribution in the reduced space. Given the parameter $\rho$, the DDPM denoises this noisy trajectory, while imposing the reduced initial condition, resulting in a prediction of the trajectory in the reduced space. This prediction can be decoded to yield a prediction in the full space, using a decoding step involving an autoencoder and POD.}
    \label{fig:1}
\end{figure}
In this work, we propose a DDPM-based ROM (DDP-ROM) that can be used to forecast the evolution over time of high-dimensional and parameterized dynamical systems efficiently. 
In particular, the DDP-ROM
\begin{itemize}
    \item is a computationally efficient ROM due to the embedding of the denoising process in a low-dimensional \emph{latent space}, as done by latent DDPMs \cite{rombach2022high}, which is built by combining linear reduction, by means of proper-orthogonal decomposition, and nonlinear autoencoders,
    \item does not rely on autoregression for predicting the evolution of the systems over time, but generates \emph{spatio-temporal} solutions directly as in \cite{ajay2022conditional, Janner2022PlanningSynthesis}, which we empirically show results in stable predictions even when extrapolating in time, outperforming a POD DL-ROM \cite{fresca2022pod} baseline and,
    \item can be used for estimating uncertainties over the estimated solutions by exploiting the probabilistic capabilities of DDPMs.
\end{itemize}
With reference to Figure \ref{fig:1}, we indicate with $\bm{y}(t_n;\bm{\rho})\in\mathbb{R}^{N_{\bm{y}}},\, n=0,\dots,H$ the high-dimensional FOM states after temporal discretization, dependent on the parameter $\boldsymbol{\rho}\in\mathcal{P}$, at different time instances $t_n\in\mathbb{R}_{\geq 0}$. We refer to the discretized states as snapshots. The reduced representation of the initial condition $\tilde{\bm{y}}(t_0;\bm{\rho})\in\mathbb{R}^{n_{\bm{y}}}$, where $n_{\bm{y}}\ll N_{\bm{y}}$, is obtained by projecting the FOM snapshot onto a low-dimensional latent space using an encoder-decoder scheme.
To generate spatio-temporal solutions in the latent space, we rely on a latent DDPM that denoises samples from a Gaussian distribution $\mathcal{N}(\bm{0}, I)$ to generate a latent space trajectory $\tilde{\bm{y}}(t_n;\boldsymbol{\rho})\in\mathbb{R}^{n_{\bm{y}}},\, n=1,\dots,H$. During generation of the latent trajectories, the parameter $\bm{\rho}$ is taken into account using conditional sampling from the DDPM, and the initial state of such trajectory is enforced by fixing $\tilde{\bm{y}}(t_0;\boldsymbol{\rho})$ at each diffusion step $k \in \{1, \dots, K\}$. After the generation of the latent space trajectory - that can have an arbitrary time horizon - the decoder maps the latent snapshots back to the high-dimensional space, resulting in a prediction of the FOM trajectory, denoted by $\hat{\bm{y}}(t_n;\boldsymbol{\rho})\in\mathbb{R}^{n_{\bm{y}}},\, n=1,\dots,H$. We test our DDP-ROM on a parameterized fluid flow around an object, even extrapolating to time instances beyond the training horizon, outperforming a well-established baseline set by deep learning-based ROMs \cite{fresca2022pod}.

The remainder of the paper is organized as follows: Section \ref{sec:preliminaries} introduces the problem setting, dimensionality reduction using POD and autoencoders, and DDPMs. Section \ref{sec:methodology} presents the proposed framework of spatio-temporal latent DDP-ROM. Section \ref{sec:numerical_experiments} shows the numerical experiments, results, and discusses the findings. We conclude in Section \ref{sec:conclusion} with an outlook on future work.

\section{Preliminaries}\label{sec:preliminaries}

Section \ref{subsec:problem_setting} starts with a description of the problem setting considered in this paper. Section \ref{subsec:Dim_Reduction} recalls important elements of dimensionality reduction, in particular POD and autoencoders \cite{hinton1993autoencoders, goodfellow2016deep}. Section \ref{subsec:DDPMs} introduces Denoising Diffusion Probabilistic Models \cite{Sohl-Dickstein2015DeepThermodynamics, Ho2020DenoisingModels}, since these are an integral part of DDP-ROM.

\subsection{Problem Setting}\label{subsec:problem_setting}
In this paper, we consider parametric initial value problems as our (high-fidelity) FOM, that arise from the spatial discretization of  nonlinear parameterized PDEs. This FOM can then be formulated by: for a given parameter vector $\boldsymbol{\rho} \in \mathcal{P}\subset \mathbb{R}^{N_{\mathcal{P}}}$ and initial condition $\bm{y}_0(\boldsymbol{\rho})$, we seek a solution $\bm{y}(\cdot;\boldsymbol{\rho}) \in C^1\left(\mathcal{I};\mathbb{R}^{N_{\bm{y}}}\right)$ such that
\begin{equation}
    \frac{d}{dt}\bm{y}(t;\boldsymbol{\rho}) = \bm{f}(t, \bm{y}(t;\boldsymbol{\rho});\boldsymbol{\rho}),  \ \ \ \ \ \bm{y}(t_0;\boldsymbol{\rho}) = \bm{y}_0(\boldsymbol{\rho}) \in \mathbb{R}^{N_{\bm{y}}},
    \label{eq:FOM}
\end{equation}
where $t \in \mathcal{I}:=(t_0, t_f]$ is the time variable and $\bm{f}:\mathcal{I} \times \mathbb{R}^{N_{\bm{y}}} \times \mathcal{P} \rightarrow \mathbb{R}^{N_{\bm{y}}}$ is a nonlinear function, encoding the spatially discretized dynamics of the nonlinear PDE.

Numerically solving the high-dimensional ODE \eqref{eq:FOM} for $N_{\boldsymbol{\rho}}$ different instances of the parameter $\boldsymbol{\rho} \in \mathcal{P}$ at various time instances yields the solution snapshots $\bm{y}(t_i,\boldsymbol{\rho}_j) \in \R^{N_{\bm{y}}}$ for $i=0,\dots, N_t-1$ and $j=1,\dots,N_{\boldsymbol{\rho}}$. Concatenating these snapshots results in a snapshot matrix $X  \in \mathbb{R}^{N_{\bm{y}} \times (N_t \times N_{\boldsymbol{\rho}})}$ with
\begin{equation}
  X =\left[ \begin{array}{ccccc} | & | & \cdots & | & | \\
         {\bm{y}(t_0;\boldsymbol{\rho}_1)} & {\bm{y}(t_1;\boldsymbol{\rho}_1)} & \cdots & {\bm{y}(t_{N_{t-2}};\boldsymbol{\rho}_{N_{\boldsymbol{\rho}}})}  &  {\bm{y}(t_{N_{t-1}};\boldsymbol{\rho}_{N_{\boldsymbol{\rho}}})} \\
         | & | & \cdots & | & |   \end{array} \right] \, .
\label{eq:data_matrix}
\end{equation}

\subsection{POD-AE Dimensionality Reduction} \label{subsec:Dim_Reduction}
In many complex scientific and engineering applications—such as fluid dynamics—the dimension of the FOM can be so enormous that simulating the system becomes prohibitively computationally expensive. In reduced-order modeling, this is mitigated by searching for a low-dimensional approximation of the state and simulating the dynamics in this reduced space instead. The dimensionality reduction method in this work is similar to the one used for POD-DL-ROM \cite{fresca2022pod} and features a two-stage approach to dimensionality reduction. The first stage is the projection onto a reduced linear basis found through proper orthogonal decomposition (POD)\cite{Sirovich1987Turbulence} and the second stage makes use of an autoencoder (AE) \cite{hinton1993autoencoders, goodfellow2016deep}. This section briefly introduces POD and AE as tools for dimensionality reduction.

POD aims to find an optimal orthogonal basis that spans a reduced space to project the solution on. In practice, this basis is often found from data by applying the singular value decomposition (SVD) to the snapshot matrix and truncating it to retain only the $r$ largest singular values and left singular vectors, yielding
\begin{equation*}
    X \approx U \Sigma V^{\top}\, ,
    \label{eq:SVD}
\end{equation*}
where $U \in \mathbb{R}^{N_{\bm{y}}\times r}$, $\Sigma \in \mathbb{R}^{r \times r}$, $V \in \mathbb{R}^{(N_t \times N_{\rho})\times r}$,  with $r$ denoting the number of retained modes. The columns of $U$ are the orthogonal basis vectors $u_i\in\mathbb{R}^{N_{\bm{y}}},i=1,\dots,r$ of the reduced space. This basis is optimal in a least-squares sense with respect to the data in the snapshot matrix. Consequently, this method for finding a reduced space is used abundantly throughout intrusive and non-intrusive MOR alike. Whereas intrusive methods would use this basis to project the FOM dynamics onto the reduced space, non-intrusive methods typically project the snapshots down using 
\begin{equation}
    \bar{\bm{y}}(t_n;\boldsymbol{\rho}_p) = U^{\top}\bm{y}(t_n;\boldsymbol{\rho}_p)\in\mathbb{R}^{r},\quad n=0,\dots,N_{t-1},\,p=1,\dots,N_{\bm{\rho}}.
\end{equation}
An approximate reconstruction of the snapshot can subsequently be recovered through 
\begin{equation*}
    \bm{y}(t_n;\boldsymbol{\rho}_p) \approx U\bar{\bm{y}}(t_n;\boldsymbol{\rho}_p),\quad n=0,\dots,N_{t-1},\,p=1,\dots,N_{\bm{\rho}}. 
\end{equation*}
The reduced snapshots can then be used to learn the reduced dynamics.

In general, the approximation error for an $N$-dimensional linear approximation space is bound by the Kolmogorov $N$-width. For advection-dominated problems, which are common in, for example, fluid dynamics, the Kolmogorov $N$-width may decay slowly with increasing $N$. Nonlinear dimensionality reduction methods can provide more efficient representations. For this reason, neural networks are a popular choice for learning nonlinear reduction mappings, since they are universal function approximators. Specifically, AEs provide a nonlinear alternative for dimensionality reduction. An AE consists of two components: (i) an encoder $\enc:\mathbb{R}^{N_{\bm{y}}} \rightarrow \mathbb{R}^{n_{\bm{y}}}$, which maps high-dimensional inputs to low-dimensional latent variables, and (ii) a decoder $\dec:\mathbb{R}^{n_{\bm{y}}} \rightarrow \mathbb{R}^{N_{\bm{y}}}$, which reconstructs the original inputs from the latent variables. Given a snapshot $\bm{y}(t_n;\boldsymbol{\rho}_p)$, we can write the AE reconstruction as
\begin{equation}
    \bm{y}(t_n;\boldsymbol{\rho}_p) \approx \dec(\enc(\bm{y}(t_n;\boldsymbol{\rho}_p);\boldsymbol{\theta}_{\enc});\boldsymbol{\theta}_{\dec})\, ,
\end{equation}
where encoder and decoder are parameterized by neural networks with learnable weights and biases indicated by $\boldsymbol{\theta}_{\enc}$ and $\boldsymbol{\theta}_{\dec}$. The learnable parameters of the AE are jointly optimized by minimizing the reconstruction error over a training dataset $Z = \{\bm{y}^{(1)}, \dots, \bm{y}^{(N)} \}$
\begin{equation}
\begin{split}
      \mathcal{L}(\boldsymbol{\theta}_{\enc}, \boldsymbol{\theta}_{\dec}) &= \frac{1}{N}\sum_{i=1}^N\left\|\bm{y}^{(i)} - \hat{\bm{y}}^{(i)}\right\|^2_2  
      = \frac{1}{N}\sum_{i=1}^N \left\|\bm{y}^{(i)} - \dec(\enc(\bm{y}^{(i)}; \boldsymbol{\theta}_{\enc}); \boldsymbol{\theta}_{\dec}) \right\|^2_2\, ,
\end{split}
\label{eq:recon_loss}
\end{equation}
where $\bm{y}^{(\text{i})}$ indicates the $i$-th sample of training data, corresponding to a FOM snapshot $\bm{y}(t_n;\boldsymbol{\rho}_p)$, and $\hat{\bm{y}}^{(\text{i})}$ its reconstruction.

For very high-dimensional systems, the training stage for AEs can be computationally expensive and data-intensive, due to the large amount of network parameters to tune. Instead \cite{fresca2022pod} proposes to first apply POD to reduce the dimensionality to an intermediate level and then train an AE on the POD coefficients. This reduces the training costs on the AE, while enabling a more efficient reduction than POD alone can provide. The autoencoder is then trained on a reduced dataset $\tilde{Z} = \{\bar{\bm{y}}^{(1)}, \dots, \bar{\bm{y}}^{(N)} \}$, analogously to the direct AE approach. In this case, the reduced representations of the snapshots are given by
\begin{equation}\label{eq:POD-AE}
    \tilde{\bm{y}}(t;\bm{\rho}) = \enc (U^\top\bm{y}(t;\bm{\rho})) \in \mathbb{R}^{n_{\bm{y}}},
\end{equation}
and the corresponding approximate reconstruction by 
\begin{equation}\label{eq:POD-AE-decode}
    \hat{\bm{y}}(t;\bm{\rho}) = U \dec (\tilde{\bm{y}}(t;\bm{\rho})) \in \mathbb{R}^{N_{\bm{y}}},
\end{equation}

\subsection{Denoising Diffusion Probabilistic Models}\label{subsec:DDPMs}
In non-intrusive MOR, the dynamics in the reduced space need to be recovered from data. In this work, we leverage the developments in generative deep-learning to train a probabilistic model that predicts the temporal evolution of an initial condition under given parameters. This is achieved through the use of DDPMs\cite{Sohl-Dickstein2015DeepThermodynamics, Ho2020DenoisingModels}. Section \ref{sec:methodology} describes the application of DDPMs to model trajectories, similar to \cite{ajay2022conditional,Janner2022PlanningSynthesis}. In this section, we briefly introduce DDPMs for general datasets.

DDPMs are a class of probabilistic generative models that learn data distributions from which new samples can be generated. DDPMs belong to the broader family of deep generative methods, alongside generative adversarial networks \cite{Goodfellow2020GenerativeNetworks}, variational autoencoders \cite{Kingma2014Auto-encodingBayes}, and score-based generative models \cite{Song2019GenerativeDistribution, Song2020ImprovedModels}. DDPMs define a parameterized memoryless stochastic process, known as a Markov chain, that progressively transforms random noise into structured data through iterative denoising.

This parameterized Markov chain is achieved as follows. DDPMs are characterized by a \emph{forward} and a \emph{reverse} diffusion process. The forward process iteratively adds Gaussian noise to a given data sample $\bm{x}\in \mathbb{R}^{N_{\bm{x}}}$ from a dataset $\mathcal{X}=\{\bm{x}_i\}_{i=1}^n$ from an unknown data distribution $p_{\text{data}}$. The result of applying the forward process to $\bm{x}$ is a sequence of samples with increasingly more noise added to them. The joint conditional probability density of this sequence given a noise-free data sample can be expressed as
\begin{equation}
q_{1:K}(\bm{x}^{1:K}|\bm{x}^0) = \prod_{k=1}^Kq_k(\bm{x}^k|\bm{x}^{k-1})\, , \ \ \ \ \ \ \  q_k(\bm{x}^k|\bm{x}^{k-1})=\mathcal{N}(\bm{x}^k; \sqrt{1-\beta_k}\bm{x}^{k-1}, \beta_k\mathbf{I}) \, ,
    \label{eq:forward_process}
\end{equation}
where the superscript $k=1, \dots, K$ is the index of the diffusion step, $K$ is the total number of diffusion steps, $\bm{x}^0=\bm{x}$ is the noise-free data, $\bm{x}^{1:K}$ is a shorthand notation for the sequence $(\bm{x}^1,\bm{x}^{2},\dots,\bm{x}^K)$, $q_k(\bm{x}^k|\bm{x}^{k-1})$ is the $k^\text{th}$ forward transition density, and $\mathcal{N}(\bm{x};\bm{\mu}, \bm{\Sigma})$ denotes a multivariate Gaussian density in $\mathbb{R}^{N_{\bm{x}}}$ with mean $\bm{\mu}\in\mathbb{R}^{N_{\bm{x}}}$ and covariance $\bm{\Sigma}\in\mathcal{S}^{N_{\bm{x}}}_+$ evaluated at $\bm{x}$. The sequence $\{\beta_k\}_{k=1}^K$, with $\beta_k\in (0, 1),k=1,\dots K$ controls the noise variance schedule of the forward diffusion process. This schedule is chosen such that after $K$ steps of the forward process, the original data distribution approaches a multivariate standard normal distribution, i.e. $\beta_k$ approaches $1$. In this way, an arbitrary data distribution $\mathcal{X}$ is approximately transformed into a distribution that is known a priori, namely a Gaussian distribution, in $K$ steps.

The reverse process aims to invert the forward process by recovering structured data from noisy samples. That is, given samples of Gaussian noise, the reverse process should generate samples according to $p_{\text{data}}$. The reverse process is learned from data using neural networks and its joint probability density is defined as
\begin{equation}
p_{\boldsymbol{\theta},0:K}(\bm{x}^{0:K}) = p_K(\bm{x}^K)\prod_{k=1}^Kp_{\bm{\theta},k}(\bm{x}^{k-1}|\bm{x}^k), \ \ \ \ \ \ \ p_{\bm{\theta},k}(\bm{x}^{k-1}|\bm{x}^k)=\mathcal{N}(\bm{x}^{k-1}; \boldsymbol{\mu}(\bm{x}^k, k; \boldsymbol{\theta}_{\boldsymbol{\mu}}), \boldsymbol{\Sigma}^k)\, ,
    \label{eq:reverse_process}
\end{equation}
where $p_K(\bm{x}^K)$ is the probability density according to which noise is sampled, i.e. a standard Gaussian, $p_{\bm{\theta},k}(\bm{x}^{k-1}|\bm{x}^k)$ is the reverse transition density, $\boldsymbol{\mu}(\bm{x}^k, k; \boldsymbol{\theta}_{\boldsymbol{\mu}})\in\mathbb{R}^{N_{\bm{x}}}$ is the mean of the reverse process at diffusion step $k$ learned by a neural network with parameters $\boldsymbol{\theta}_{\boldsymbol{\mu}}\in\mathbb{R}^{N_{\bm{\theta}}}$, and $\boldsymbol{\Sigma}^k\in\mathcal{S}^{N_{\bm{x}}}_{++}$ is a predefined covariance (often set via a cosine schedule \cite{Nichol2021ImprovedModels}).

\subsubsection{Training Objectives} 
DDPMs can be trained by maximizing a variational lower bound, similarly to variational autoencoders \cite{Sohl-Dickstein2015DeepThermodynamics}. A widely used alternative is the simplified noise-prediction objective introduced in \cite{Ho2020DenoisingModels}, where the model learns to estimate the noise injected during diffusion rather than directly predict the mean of the reverse process from \eqref{eq:reverse_process}. This can be done by making use of the closed-form of the forward process
\begin{equation}
    \bm{x}^k = \sqrt{\bar{\alpha}_k}\bm{x}^0+\sqrt{1-\bar{\alpha}_k}\boldsymbol{\epsilon},
\end{equation}
where $\boldsymbol{\epsilon} \sim \mathcal{N}(\bm{0}, I)$ is the injected noise, $\Bar{\alpha}_k=\prod_{j=1}^k\alpha_j$ with $\alpha_k=1-\beta_k$. The training objective is then
\begin{equation}
\mathcal{L}(\boldsymbol{\theta}_{\boldsymbol{\epsilon}}) = \mathbb{E}_{k, \boldsymbol{\epsilon}, \bm{x}^0} \big[\left\|\boldsymbol{\epsilon} - \tilde{\bm{\epsilon}}(\bm{x}^k, k; \boldsymbol{\theta}_{\boldsymbol{\epsilon}}) \right\|_2^2\big]\, ,
    \label{eq:diffusion_simple_objective}
\end{equation}
where $\boldsymbol{\theta}_{\boldsymbol{\epsilon}}\in\mathbb{R}^{N_{\bm{\theta}}}$ indicates the learnable parameters of the noise-prediction neural network $\tilde{\bm{\epsilon}}$. Using the closed-form forward process, the loss can be written as
\begin{equation}
    \mathbb{E}_{k, \boldsymbol{\epsilon}, \bm{x}^0} \big[\left\|\boldsymbol{\epsilon} - \tilde{\boldsymbol{\epsilon}}(\sqrt{\bar{\alpha_k}}\bm{x}^0+\sqrt{1-\bar{\alpha}_k}\boldsymbol{\epsilon}, k;\boldsymbol{\theta}_{\boldsymbol{\epsilon}}) \right\|_2^2\big]\, ,
\end{equation}
where $k \sim \mathcal{U}\{1, \dots, K \}$ with $\mathcal{U}$ indicating a uniform distribution. 
Under this formulation, the mean of the reverse process can be rewritten as:
\begin{equation}
\boldsymbol{\mu}(\bm{x}^k, k)=\frac{1}{\sqrt{\alpha}_k}\Big(\bm{x}^k-\frac{1-\alpha_k}{\sqrt{1-\bar{\alpha}_k}}\tilde{\boldsymbol{\epsilon}}(\bm{x}^k, k; \boldsymbol{\theta}_{\boldsymbol{\epsilon}})\Big).
    \label{eq:mu_rewritten}
\end{equation}

\subsubsection{Conditional Generation}\label{subsec:conditional_generation}
After training, DDPMs can generate samples from the training distribution by sampling from a standard normal distribution and applying the learned reverse process to iteratively denoise the samples. However, directly applying the reverse process yields a random sample matching the training data distribution, without regard for the characteristics that one would like the sample to have. For instance, in image generation one may want to generate an image of a cat, rather than any image at all, and for predicting dynamical systems one may want to limit sample generation to certain dynamical parameters. In DDPMs, there are two main strategies for imposing such conditions on the generative process, namely (i) the \emph{classifier guidance} \cite{Dhariwal2021DiffusionSynthesis} and the \emph{classifier-free guidance} \cite{ho2022classifier}. 

Classifier guidance conditions the reverse process using the gradients of an auxiliary classifier trained to predict labels from noisy samples. While effective, this approach requires training an additional model. Conversely, classifier-free guidance avoids auxiliary classifiers by incorporating conditioning variables directly into the diffusion model. Let $\bm{c} \in \mathbb{R}^{N_{\bm{c}}}$ denote the conditioning input, the guided noise prediction is computed as
\begin{equation}
    \hat{\boldsymbol{\epsilon}} = \tilde{\boldsymbol{\epsilon}}(\bm{x}^k, \varnothing, k; \boldsymbol{\theta}_{\tilde{\boldsymbol{\epsilon}}}) + \omega(\tilde{\boldsymbol{\epsilon}}(\bm{x}^k, \bm{c}, k; \boldsymbol{\theta}_{\boldsymbol{\epsilon}}) - \tilde{\boldsymbol{\epsilon}}(\bm{x}^k, \varnothing, k; \boldsymbol{\theta}_{\boldsymbol{\epsilon}}))\, ,
    \label{background:classifier_free:combined_score}
\end{equation}
where $\omega \in \R$ controls the trade-off between unconditional and conditional generation, and $\varnothing$ denotes a null-token that is passed to the noise prediction network when no conditional information is to be used for the prediction\footnote{The null-token $\varnothing$ should not be interpreted as the zero element in $\mathbb{R}^{N_{\bm{c}}}$, as this would refer to conditioning on $\bm{c}=0$. In practice, it is implemented by learning an embedding of $\bm{c}$ and choosing $\varnothing$ to correspond to $\bm{0}$ in this embedding space. The embedding is then passed as an input to the noise prediction network.}. Classifier-free training modifies the objective in \eqref{eq:diffusion_simple_objective} by randomly dropping the conditioning input: 
\begin{equation}
\begin{split}
\mathcal{L}(\boldsymbol{\theta}_{\boldsymbol{\epsilon}}) &= \mathbb{E}_{k, \boldsymbol{\epsilon}, \bm{x}^0, \eta \sim \text{Bern}(\pi)} = \big[\left\|\boldsymbol{\epsilon} - \tilde{\boldsymbol{\epsilon}}(\boldsymbol{\tau}^k, (1-\eta)\bm{c} + \eta \varnothing, k; \boldsymbol{\theta}_{\boldsymbol{\epsilon}}) \right\|^2\big]\, , \\
&= \mathbb{E}_{k, \boldsymbol{\epsilon}, \boldsymbol{\tau}^0, \eta \sim \text{Bern}(\pi)} = \big[\left\|\boldsymbol{\epsilon} - \tilde{\boldsymbol{\epsilon}}(\sqrt{\bar{\alpha}_k}\boldsymbol{\tau}^0+\sqrt{1-\bar{\alpha}_k}\boldsymbol{\epsilon}, (1-\eta)\bm{c} + \eta \varnothing, k; \boldsymbol{\theta}_{\boldsymbol{\epsilon}}) \right\|^2\big]\, ,
\end{split}
    \label{eq:classifier-free_diffusion_simple_objective}
\end{equation}
where $\eta \sim \text{Bern}(\pi)$ determines whether the conditioning is used ($\eta=0$) or dropped ($\eta=1$).

Another important conditioning mechanism is \emph{inpainting}, where parts of the data are known and fixed while the model reconstructs the missing components. During denoising, the reverse diffusion process is constrained to preserve observed regions and generate consistent completions in unobserved areas. Producing coherent interpolations under such constraints is challenging, making inpainting one of the most extensively studied problems in diffusion modeling \cite{yang2023diffusion}.

\section{Spatio-temporal Latent Denoising Diffusion Reduced-order Model}\label{sec:methodology}

We propose a novel framework for constructing ROMs using a latent DDPM \cite{rombach2022high} which we call DDP-ROM. DDP-ROM \emph{(i)} operates in a low-dimensional latent space, obtained with linear or nonlinear dimensionality-reduction techniques (see Section \ref{subsec:Dim_Reduction}), to reduce the computational cost of the reverse diffusion process for high-dimensional data, and  \emph{(ii)} generates spatio-temporal solutions conditioned on an initial state and an instance of the parameter vector $\boldsymbol{\rho}$ to mitigate the cost of performing autoregressive rollouts and improve temporal consistency of the generated trajectories. 

\subsection{Offline and Online Phases}
Given the training dataset composed by the FOM snapshots from \eqref{eq:data_matrix}, we define a trajectory $\boldsymbol{\tau}_p$ as
\begin{equation}
    \boldsymbol{\tau}_p = [\bm{y}(t_0;\boldsymbol{\rho}_p), \dots, \bm{y}(t_n;\boldsymbol{\rho}_p), \dots, \bm{y}(t_{N_t-1};\boldsymbol{\rho}_p)]\in\mathbb{R}^{N_{\bm{y}}\times N_{t}},
    \label{eq:trajectory}
\end{equation}
with $p\in\{1,\dots,N_{\bm{\rho}}\}$.
The dimensionality of each snapshot $\bm{y}(t_n;\boldsymbol{\rho}_p)$ and consequently of the trajectory $\boldsymbol{\tau}_p$ can be large if the data describe complex systems with high degree of accuracy, e.g., the FOM uses a fine mesh. Evaluating a model that predicts these trajectories directly in $\mathbb{R}^{N_y}$ may be prohibitively computationally expensive. Therefore, we construct a reduced space of dimension $n_{\bm{y}} \ll N_y$ in a two-stage approach, as introduced in Section \ref{subsec:Dim_Reduction}. Firstly, we construct a POD-basis from the snapshot matrix (from \eqref{eq:data_matrix}). Secondly, to provide additional nonlinear reduction, we project the snapshots onto the POD-basis and use the POD coefficients to train an AE (using the reconstruction loss in \eqref{eq:recon_loss}). The reduced representation of a snapshot is found through \eqref{eq:POD-AE}.

A latent trajectory is now composed of these projected snapshots:
\begin{equation}
    \tilde{\boldsymbol{\tau}}_p = [\tilde{\bm{y}}(t_0;\boldsymbol{\rho}_p), \cdots, \tilde{\bm{y}}(t_n;\boldsymbol{\rho}_p), \cdots, \tilde{\bm{y}}(t_{N_t-1};\boldsymbol{\rho}_p)]\in\mathbb{R}^{n_{\bm{y}}\times N_{t}}.
    \label{eq:latent_trajectory}
\end{equation}
Then, to model the dynamics in the reduced space, the collection of latent trajectories, rather than the collection of individual snapshots, is used to train a DDPM. In particular, we can rewrite the forward process (see \eqref{eq:forward_process}) as:
\begin{equation}
q_{1:K}(\tilde{\boldsymbol{\tau}}_p^{1:K}|\tilde{\boldsymbol{\tau}}_p^0) = \prod_{k=1}^K q_k(\tilde{\boldsymbol{\tau}}_p^k|\tilde{\boldsymbol{\tau}}_p^{k-1})\, , \ \ \ \ \  q_k(\tilde{\boldsymbol{\tau}}_p^k|\tilde{\boldsymbol{\tau}}^{k-1})=\mathcal{N}(\tilde{\boldsymbol{\tau}}_p^k; \sqrt{1-\beta_k}\tilde{\boldsymbol{\tau}}_p^{k-1}, \beta_k\mathbf{I}) \, ,
    \label{eq:forward_process_tau}
\end{equation}
and the neural network-parameterized reverse process (from \eqref{eq:reverse_process}) as:
\begin{equation}
p_{\boldsymbol{\theta},0:K}(\tilde{\boldsymbol{\tau}}_p^{0:K}) = \prod_{k=1}^K p_{\boldsymbol{\theta},k}(\tilde{\boldsymbol{\tau}}_p^{k-1}|\tilde{\boldsymbol{\tau}}_p^k), \ \ \ \ \ p_{\boldsymbol{\theta},k}(\tilde{\boldsymbol{\tau}}_p^{k-1}|\tilde{\boldsymbol{\tau}}_p^k)=\mathcal{N}(\tilde{\boldsymbol{\tau}}_p^{k-1}; \boldsymbol{\mu}(\tilde{\boldsymbol{\tau}}_p^k, k), \boldsymbol{\Sigma}^k)\, ,
    \label{eq:reverse_process_tau}
\end{equation}
where $\tilde{\boldsymbol{\tau}}_p^k$ indicate the $p^{\text{th}}$ trajectory at the $k^{\text{th}}$ diffusion step. To be able to predict trajectories of given parameters, we use a classifier-free approach (see Section \ref{subsec:conditional_generation}), such that the DDP-ROM learns the conditional relation between the parameters and latent trajectories. We parametrize the reverse process using a noise-prediction neural network with learnable parameters indicated by $\boldsymbol{\theta}_{\boldsymbol{\epsilon}}$
\begin{equation}
\begin{split}
    \boldsymbol{\mu}(\tilde{\boldsymbol{\tau}}^k_p, k) &= \frac{1}{\sqrt{\alpha}_k}\Big(\tilde{\boldsymbol{\tau}}^k_p-\frac{1-\alpha_k}{\sqrt{1-\bar{\alpha}_k}}\tilde{\boldsymbol{\epsilon}}(\tilde{\boldsymbol{\tau}}^k_p, k; \boldsymbol{\theta}_{\boldsymbol{\epsilon}})\Big)\, , \\
    \tilde{\boldsymbol{\epsilon}}(\tilde{\boldsymbol{\tau}}^k_p, k; \boldsymbol{\theta}_{\boldsymbol{\epsilon}}) &= \tilde{\boldsymbol{\epsilon}}(\tilde{\boldsymbol{\tau}}_p^k, \bm{c}=\varnothing, k; \boldsymbol{\theta}_{\boldsymbol{\epsilon}}) + \omega(\tilde{\boldsymbol{\epsilon}}(\tilde{\boldsymbol{\tau}}_p^k, \bm{c}=\boldsymbol{\rho}_p, k; \boldsymbol{\theta}_{\boldsymbol{\epsilon}}) - \tilde{\boldsymbol{\epsilon}}(\tilde{\boldsymbol{\tau}}_p^k, \bm{c}=\varnothing, k; \boldsymbol{\theta}_{\boldsymbol{\epsilon}}))\, , \\ 
\end{split}
\label{eq:cond_sampling_tau}
\end{equation}
where the conditioning input $\bm{c}$ corresponds to the parameters $\boldsymbol{\rho}_p$.
We train DDP-ROM using the classier-free objective in \eqref{eq:classifier-free_diffusion_simple_objective} that we rewrite as:
\begin{equation}
\begin{split}
\mathcal{L}(\boldsymbol{\theta}_{\boldsymbol{\epsilon}}) &= \mathbb{E}_{k, \boldsymbol{\epsilon}, \tilde{\boldsymbol{\tau}}^0_p, \eta \sim \text{Bern}(\pi)} = \big[\left\|\boldsymbol{\epsilon} - \boldsymbol{\epsilon}(\tilde{\boldsymbol{\tau}}^k_p, (1-\eta)\bm{c} + \eta \varnothing, k; \boldsymbol{\theta}_{\boldsymbol{\epsilon}}) \right\|^2\big]\, , \\
&= \mathbb{E}_{k, \boldsymbol{\epsilon}, \tilde{\boldsymbol{\tau}}^0_m, \eta \sim \text{Bern}(\pi)} = \big[\left\|\boldsymbol{\epsilon} - \boldsymbol{\epsilon}(\sqrt{\bar{\alpha_k}}\tilde{\boldsymbol{\tau}}^0_p+\sqrt{1-\bar{\alpha}_k}\boldsymbol{\epsilon}, (1-\eta)\bm{c} + \eta \varnothing, k; \boldsymbol{\theta}_{\boldsymbol{\epsilon}}) \right\|^2\big].
\end{split}
    \label{eq:classifier-free_DDPMPDE_objective}
\end{equation}

In the online phase, we can sample a random noise vector $\check{\boldsymbol{\tau}}^K_H\sim \mathcal{N}(\mathbf{0}, \mathbf{I}) \in \mathbb{R}^{n_{\bm{y}}\times H}$, with $H$ equal to the prediction horizon and we can utilize the reverse diffusion process to generate spatio-temporal latent solutions
\begin{equation}
    \check{\boldsymbol{\tau}}_H(\bm{\rho}) =  [\check{\bm{y}}(t_0;\boldsymbol{\rho}), \cdots, \check{\bm{y}}(t_n;\boldsymbol{\rho}), \cdots, \check{\bm{y}}(t_{H-1};\boldsymbol{\rho})]\in\mathbb{R}^{n_{\bm{y}}\times H},
\end{equation}
 conditioned on (new) instances of the parameter $\boldsymbol{\rho}\in\mathcal{P}$. The latent solutions can then be projected back to the original space using the decoder and POD basis using \eqref{eq:POD-AE-decode} to obtain the predicted snapshots $\hat{\bm{y}}(t_i;\bm{\rho}),i=0,\dots,H-1$ and the predicted trajectory $\hat{\bm{\tau}}(\bm{\rho})$. It is worth highlighting the $H$ does not need to coincide with $N_t$. In fact, in our numerical experiments, we will show that by choosing $\boldsymbol{\epsilon}(\boldsymbol{\tau}, k; \boldsymbol{\theta}_{\boldsymbol{\epsilon}})$ to be convolutional along the time dimension of $\boldsymbol{\tau}$, $H$ can be larger than $N_t$ without affecting the performance of DDP-ROM, i.e., DDP-ROM is capable of extrapolating in time. In Algorithm \ref{alg:online_phase}, we show the pseudo-code of the online phase of DDP-ROM.
\begin{algorithm}
\caption{Online Phase}\label{alg:online_phase}
\begin{algorithmic}[1]
\State Given a trained noise-prediction model $\boldsymbol{\epsilon}(\cdot;\boldsymbol{\theta}_{\boldsymbol{\epsilon}})$, the reduced POD basis $U$, and a trained AE $(\enc, \dec)$.
\State Choose prediction horizon $H$ and parameter instance $\boldsymbol{\rho}$. \State Project initial state $\bm{y}(t_0;\boldsymbol{\rho})$ using $ \enc \circ U^\top$  to its latent representation $\tilde{\bm{y}}(t_0;\boldsymbol{\rho})$
\State Sample $\check{\boldsymbol{\tau}}^K_H\sim \mathcal{N}(\mathbf{0}, \mathbf{I}) \in \mathbb{R}^{n_{\bm{y}}\times H}$
\For{$k=K,\cdots,1$}
\State Fix initial state of the latent trajectory to be $\tilde{\bm{y}}(t_0;\boldsymbol{\rho})$
\State  $\hat{\boldsymbol{\epsilon}} = \boldsymbol{\epsilon}(\check{\boldsymbol{\tau}}^k_H, \bm{c}=\varnothing, k; \boldsymbol{\theta}_{\boldsymbol{\epsilon}}) + \omega(\boldsymbol{\epsilon}(\check{\boldsymbol{\tau}}^k_H, \bm{c}=\boldsymbol{\rho}, k; \boldsymbol{\theta}) - \boldsymbol{\epsilon}(\check{\boldsymbol{\tau}}^k_H, \bm{c}=\varnothing, k; \boldsymbol{\theta}_{\boldsymbol{\epsilon}}))$ \Comment{Conditional Sampling}
\State  $\boldsymbol{\mu}(\check{\boldsymbol{\tau}}^k_H, k) = \frac{1}{\sqrt{\alpha}_k}\Big(\check{\boldsymbol{\tau}}^k_H-\frac{1-\alpha_k}{\sqrt{1-\bar{\alpha}_k}}\boldsymbol{\epsilon}(\check{\boldsymbol{\tau}}^k_H, \bm{c}, k; \boldsymbol{\theta}_{\boldsymbol{\epsilon}})\Big)$
\State $\check{\boldsymbol{\tau}}^{k-1}_H \sim \mathcal{N}(\check{\boldsymbol{\tau}}^{k-1}_H; \boldsymbol{\mu}(\check{\boldsymbol{\tau}}^k_H, k), \boldsymbol{\Sigma}^k)$
\EndFor
\State Obtain $\hat{\boldsymbol{\tau}}$ from $\check{\boldsymbol{\tau}}_H$ using $U \circ \dec$
\end{algorithmic}
\end{algorithm}

\section{Numerical Experiments}\label{sec:numerical_experiments}
We test the forecasting capabilities of DDP-ROM on a 2-dimensional fluid flow around an obstacle. We compare the performance of DDP-ROM with POD DL-ROM, (see Appendix \ref{app:dl-rom} for details on POD DL-ROMs).

After training, we evaluate the prediction abilities of the models performance using various error metrics. The prediction is a relative error for a single parameter instance at a given time defined as
\begin{equation} 
        e_{\text{pred}}\left(t; \bm{\rho}\right) = \sqrt{ \frac{\left\|\bm{y}(t;\boldsymbol{\rho}) - \hat{\bm{y}}(t;\boldsymbol{\rho})\right\|^2}{\left\|\bm{y}(t;\boldsymbol{\rho})\right\|^2 } }.
        \label{eq:prediction_error}
\end{equation}
We define the mean prediction error as the mean of the prediction error over all parameter instances in the test data set
\begin{equation}
        e_{\text{mean}}\left(t\right) = \frac{1}{N_{\text{test}}} \sum_{i=1}^{N_{\text{test}}} e_{\text{pred}}\left(t;\bm{\rho}_i\right).
        \label{eq:mean_prediction_error}
\end{equation}
The total prediction error assesses the overall performance of a ROM over all $N_{\text{test}}$ trajectories in the test set and is defined as
\begin{equation}
        e_{\text{tot}} = \frac{1}{N_{\text{test}}} \sum_{i=1}^{N_{\text{test}}}\sqrt{ \frac{\sum_{k=1}^{N_t}\left\|\bm{y}(t_k;\boldsymbol{\rho}_i) - \hat{\bm{y}}(t_k;\boldsymbol{\rho}_i)\right\|^2}{\sum_{k=1}^{N_t}\left\|\bm{y}(t_k;\boldsymbol{\rho}_i)\right\|^2 } }.
        \label{eq:tot_prediction_error}
    \end{equation}
Since DDP-ROM is a probabilistic model, there may be some variability in the predictions. The mean absolute prediction error measures the error of a snapshot prediction over the domain, where the mean is with respect to the $N_{\text{repeat}}$ predictions of the same parameters. It is defined as
    \begin{equation}
        \boldsymbol{e}_{\text{abs}}\left(t;\bm{\rho}\right) = \frac{1}{N_{\text{repeat}}} \sum_1^{N_{\text{repeat}}}|\bm{y}(t;\boldsymbol{\rho}) - \hat{\bm{y}}(t;\boldsymbol{\rho})|.
    \label{eq:abs_error}
    \end{equation}
Additionally, we look at the kinetic energy defined as a function of the velocity $\bm{v}$
\begin{equation}\label{eq:e_kin}
    E_{\text{kin}}(t)=\frac{1}{2}\int_{\Omega} \boldsymbol{v}(\bm{x},t)\cdot \boldsymbol{v}(\bm{x},t)\,\mathrm{d}V,
\end{equation}
which is computed on the finite element mesh.

\subsection{Fluid Flow Around Obstacle}\label{subsec:fluidflow}
As a test case, we consider the problem of learning a ROM for a fluid flow around an obstacle. We assume incompressible flow, whose dynamics are described by the unsteady Navier-Stokes equations:
\begin{equation}
    \begin{aligned}
    \begin{cases}
        \frac{\partial\bm{v}}{\partial t}(\bm{x},t) - \nu \Delta \bm{v}(\bm{x},t) + (\bm{v}(\bm{x},t) \cdot \nabla )\bm{v}(\bm{x},t) + \nabla p(\bm{x},t) = 0 \quad &(\bm{x},t)\in\Omega\times(0,T] \\
        \nabla \cdot \bm{v}(\bm{x},t) = 0 \quad &(\bm{x},t)\in\Omega\times(0,T]\\
        \bm{v}(\bm{x},t)=\bm{v}_{\mathrm{in}}(\bm{x})  \quad &(\bm{x},t)\in\partial\Omega_{\mathrm{in}}\times(0,T]\\
        \frac{\partial\bm{v}}{\partial x_2}(\bm{x},t)=\bm{0}  \quad &(\bm{x},t)\in\partial\Omega_{\mathrm{wall}}\times(0,T]\\
        \bm{v}(\bm{x},t)=\bm{0}  \quad &(\bm{x},t)\in\partial\Omega_{\mathrm{object}}\times(0,T]\\
         p(\bm{x},t)=0  \quad &(\bm{x},t)\in\partial\Omega_{\mathrm{out}}\times(0,T]\\
        \bm{v}(\bm{x},0)=\bm{v}_0(\bm{x}) \quad &\bm{x}\in\Omega\\
    \end{cases}
    \end{aligned}
\end{equation}
where $\Omega \subset \mathbb{R}^2$ is the spatial domain, $\bm{v}:\Omega\times[0, T] \rightarrow\mathbb{R}^2$ is the fluid velocity, and $p:\Omega \times [0,T] \rightarrow \mathbb{R}$ the pressure. The boundary of the domain consists of the left side $\partial\Omega_{in}$, right side $\partial\Omega_{out}$, and the walls at the top, bottom $\partial\Omega_{wall}$, and the object boundary $\partial\Omega_{object}$, see Figure \ref{fig:domain}. 

\tikzset{
    boundary/.style={line width=1.1pt, blue!70!black},
    object/.style={line width=1.0pt, red!75!black, fill=red!12},
    label/.style={font=\small, fill=none, inner sep=2pt},
}
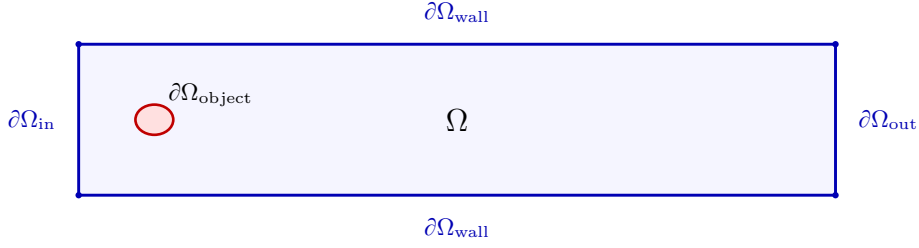
\begin{figure}
    \centering
    \begin{tikzpicture}[node distance=5mm]
    \coordinate (bottomleft) at (0,0);
    \coordinate (bottomright) at (10,0);
    \coordinate (topright) at (10,2);
    \coordinate (topleft) at (0,2);
    \coordinate (objectcenter) at (1,1);

    \fill[blue!4] (bottomleft) rectangle (topright);

    \draw[boundary] (bottomleft) -- (topleft)
        node[pos=0.5, left=6pt, label] {$\partial\Omega_{\mathrm{in}}$};
    \draw[boundary] (bottomright) -- (topright)
        node[pos=0.5, right=6pt, label] {$\partial\Omega_{\mathrm{out}}$};
    \draw[boundary] (topleft) -- (topright)
        node[pos=0.5, above=6pt, label] {$\partial\Omega_{\mathrm{wall}}$};
    \draw[boundary] (bottomleft) -- (bottomright)
        node[pos=0.5, below=6pt, label] {$\partial\Omega_{\mathrm{wall}}$};

    \draw[object] (objectcenter) ellipse (0.25cm and 0.2cm);
    \node[label, above right=3pt] at (1,1) {$\partial\Omega_{\mathrm{object}}$};

    \node[label,centered] at (5,1) {\large $\Omega$};

    \fill[blue!70!black] (bottomleft) circle (1.2pt);
    \fill[blue!70!black] (bottomright) circle (1.2pt);
    \fill[blue!70!black] (topright) circle (1.2pt);
    \fill[blue!70!black] (topleft) circle (1.2pt);
\end{tikzpicture}
    \caption{Spatial domain $\Omega$ of the test case and its boundaries $\partial\Omega_{\mathrm{in}}$, $\partial\Omega_{\mathrm{out}}$, $\partial\Omega_{\mathrm{wall}}$, and $\partial\Omega_{\mathrm{object}}$.}
    \label{fig:domain}
\end{figure}

The initial condition is given by $\bm{v}_0(\bm{x})=\bm{0}$. Fluid velocity on the left boundary has a flow velocity $\gamma_{\mathrm{in}} \in [1, 6]$ and angle of attack $\alpha_{\mathrm{in}} \in [-1, 1]$, such that
\begin{equation}
    \bm{v}_{\mathrm{in}}(\bm{x}) = 
        \big(\gamma_{\mathrm{in}}\cos{\alpha_{\mathrm{in}}}, x_2(2-x_2)\gamma_{\mathrm{in}}\sin{\alpha_{\mathrm{in}}}\big)^\top\,,
\end{equation}
where the parabolic profile $x_2(2-x_2)$ is used to prevent discontinuities. On the object is a no-slip condition, on the right boundary is a constant pressure condition and on the top and bottom walls is a free slip condition. We consider a 2-dimensional channel with dimensions $[0, 10]\times[0, 2]$ with an ellipsoidal object centered in $(1,1)$ and with radius equal to $0.25$ in $x$-direction and $0.20$ in $y$-direction. We rely on a mesh with $N=80592$ degrees of freedom in the velocity and we solve the FOM with a $\Delta t=0.05$s using FEniCS 2019.1.0 \cite{alnaes2015fenics}. The training dataset consists of $800$ trajectories with length equal to $T=10$s, where the parameter vector $\boldsymbol{\rho}=[\gamma_{\mathrm{in}}, \alpha_{\mathrm{in}}]$ is randomly sampled from uniform distributions $\gamma_{\mathrm{in}}\sim \mathcal{U}(1, 6)$ and $ \alpha_{\mathrm{in}}\sim \mathcal{U}(-1,1)$ for each trajectory. The test dataset consists of $125$ trajectories with different parameters sampled uniformly from the parameter space, but unseen during training. To assess the time-extrapolation capabilities of DDP-ROM, the test set trajectories have length equal to $T=20$s.

\begin{figure}[h!]
    \centering
    \makebox[\textwidth][c]{\includegraphics[width=0.6\textwidth]{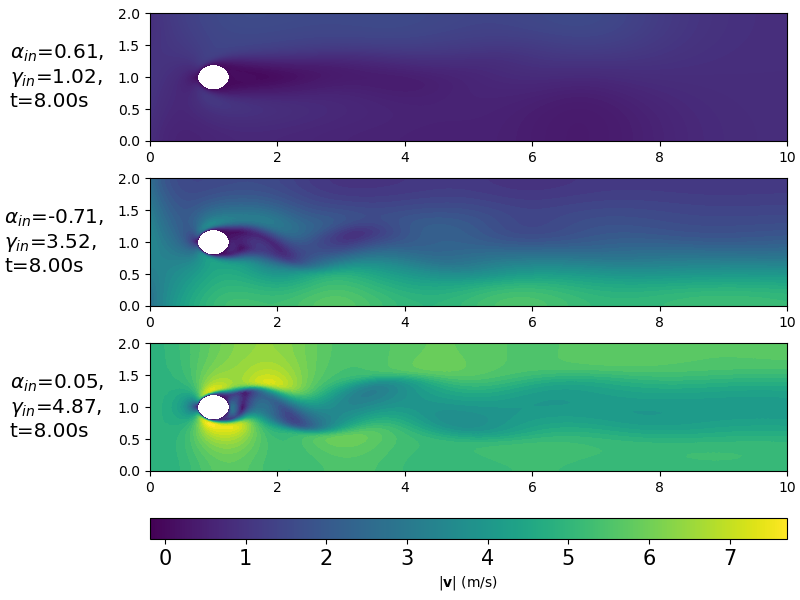}}
    \caption{Example snapshots showing velocity magnitude for three sets of parameters (indicated on the left side of the plot) at $t=8$s, showcasing the variability in the data set.}
    \label{fig:flow_examples}
\end{figure}

In Figure \ref{fig:flow_examples} we show some examples of solution snapshots of the velocity magnitudes for various parameter instances, namely $\alpha_{\mathrm{in}}, \gamma_{\mathrm{in}}$ and for $t=8$s.

Due to the high-dimensionality of the snapshots, we first project the velocity snapshots to a lower dimension using POD. A separate POD basis is used for the two velocity directions. The POD bases are computed using the \texttt{randomized\_svd} functionality from the \texttt{scikit-learn}\cite{scikit-learn} library. For each direction, the reduced dimension is $r=256$, which accounts for more than $99.95\%$ of signal energy of the training data set. The POD coordinates are concatenated to obtain a $512$-dimensional state representation. Then, we further reduce the dimensionality of the POD coefficients using an AE. We test three different AE architectures:
\begin{enumerate}
    \item mlp-AE composed of feedforward layers,
    \item conv-AE composed of convolutional layers as proposed in the original POD DL-ROM work \cite{fresca2022pod}, and
    \item inv-AE composed of invertible neural network layers as proposed in \cite{botteghi2026deep}.
\end{enumerate}
We will compare DDP-ROM and POD DL-ROM using the same POD-basis dimension $r$ and same AE latent dimensions $n_{\bm{y}} \in \{2, 3, 4, 5, 10, 20 \}$. For details about the neural networks, see Appendix \ref{app:neural_network}.

On a laptop equipped with an Intel Core Ultra 7 155H CPU and NVIDIA RTX 1000 Ada Generation Laptop GPU, the time it takes to simulate a trajectory of $400$ time steps, corresponding to $20$s has been recorded. The FOM takes $11.5$ minutes to generate one trajectory. Depending on the reduced dimension $n_{\bm{y}}$, DDP-ROM takes on average between $1.68$s and $1.75$s and DLROM between $4.34\cdot10^{-2}$s and $5.28\cdot 10^{-2}$s. DLROM is faster because of the simple architecture of the dynamical model, which causes computation time to be dominated by the POD expansion. Instead, DDP-ROM's computation time is dominated by the DDPM inference. Note, however, that there exist several methods to speed up inference of DDPMs, which have not been implemented in this work, see for example \cite{Ma2025Efficient}.

Table \ref{tab:err_lat} shows the total prediction error (see \eqref{eq:tot_prediction_error}) over the entire test set for DDP-ROM and POD DL-ROM using different autoencoder architectures and latent dimensions. Regardless of the autoencoder architecture, every DDP-ROM outperforms all POD DL-ROMs at almost all latent dimensions.

\begin{table}[t]
\centering
\setlength{\tabcolsep}{4pt}
\caption{Total prediction error computed over the test set of unseen parameters and time extrapolation for DDP-ROM and DL-ROM with feedforward AE (mlp-AE), convolutional AE (conv-AE), and invertible AE (inv-AE).}\label{tab:err_lat}
\rowcolors{2}{gray!10}{white}
\begin{tabular}{l l l l l l l}
\toprule
Lat. dim. 
  & \multicolumn{1}{c}{2}
  & \multicolumn{1}{c}{3}
  & \multicolumn{1}{c}{4}
  & \multicolumn{1}{c}{5}
  & \multicolumn{1}{c}{10}
  & \multicolumn{1}{c}{20} \\
\midrule
DDP-ROM mlp-AE & $\boldsymbol{1.29 \cdot 10^{-1}}$ & $1.56 \cdot 10^{-1}$ & $1.51 \cdot 10^{-1}$ & $1.41 \cdot 10^{-1}$ & $1.21 \cdot 10^{-1}$ & $1.48 \cdot 10^{-1}$ \\
DDP-ROM conv-AE & $1.37 \cdot 10^{-1}$ & \boldsymbol{$1.42 \cdot 10^{-1}$} & $1.54 \cdot 10^{-1}$ & $1.43 \cdot 10^{-1}$ & \boldsymbol{$1.18 \cdot 10^{-1}$} & \boldsymbol{$1.31 \cdot 10^{-1}$} \\
DDP-ROM inv-AE & $1.32 \cdot 10^{-1}$ & $1.51 \cdot 10^{-1}$ & \boldsymbol{$1.36 \cdot 10^{-1}$} & \boldsymbol{$1.26 \cdot 10^{-1}$} & $1.27 \cdot 10^{-1}$ & $1.45 \cdot 10^{-1}$ \\
POD DL-ROM mlp-AE & $5.22 \cdot 10^{-1}$ & $2.11 \cdot 10^{0}$ & $9.90 \cdot 10^{-1}$ & $7.58 \cdot 10^{-1}$ & $1.27 \cdot 10^{1}$ & $3.84 \cdot 10^{0}$ \\
POD DL-ROM conv-AE & $3.61 \cdot 10^{-1}$ & $2.07 \cdot 10^{0}$ & $1.20 \cdot 10^{1}$ & $1.43 \cdot 10^{1}$ & $5.57 \cdot 10^{0}$ & $9.65 \cdot 10^{0}$ \\
POD DL-ROM inv-AE & $2.06 \cdot 10^{-1}$ & $1.50 \cdot 10^{-1}$ & $1.95 \cdot 10^{-1}$ & $1.80 \cdot 10^{-1}$ & $6.79 \cdot 10^{-1}$ & $7.01 \cdot 10^{-1}$ \\
\bottomrule
\end{tabular}
    
\end{table}

\begin{figure}[h]
    \centering
    \includegraphics[width=\textwidth]{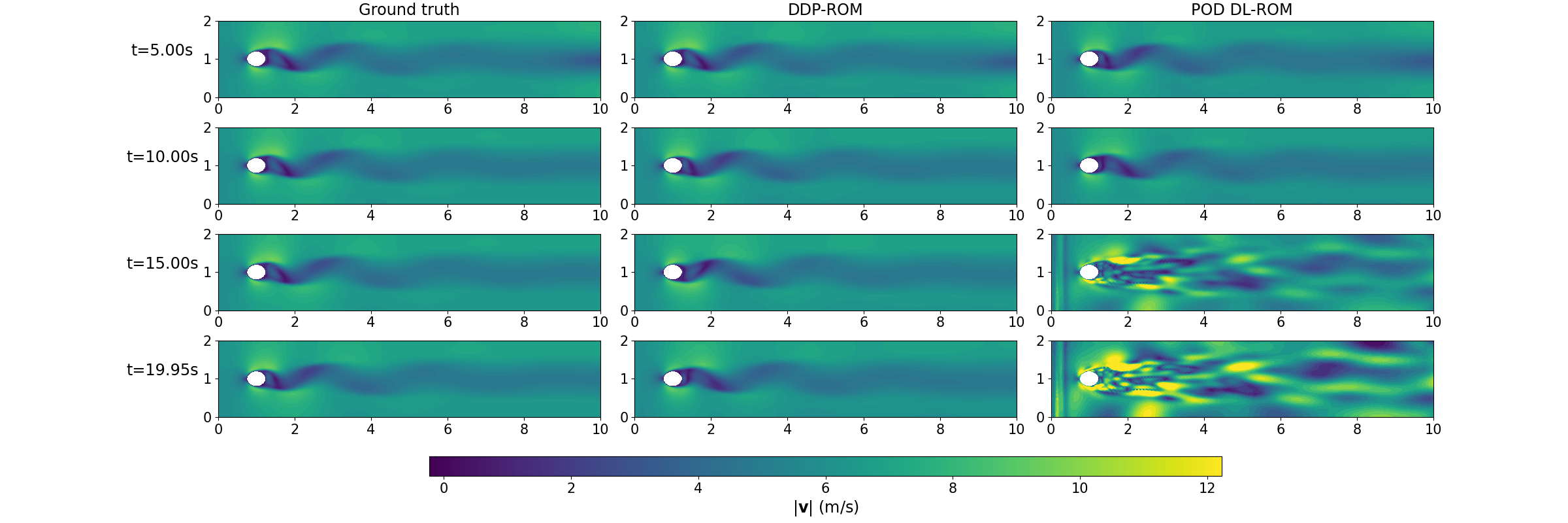}
    \caption{Velocity magnitude snapshots for a trajectory in the test set with $\alpha_{\text{in}}=0.05$ rad and $\gamma_{\text{in}}=5.99$ m/s. Left column is the ground truth, middle column the DDP-ROM prediction, and right column is DL-ROM prediction.}
    \label{fig:vmag_t11}
\end{figure}

\begin{figure}[h]
    \centering
    \includegraphics[width=\textwidth]{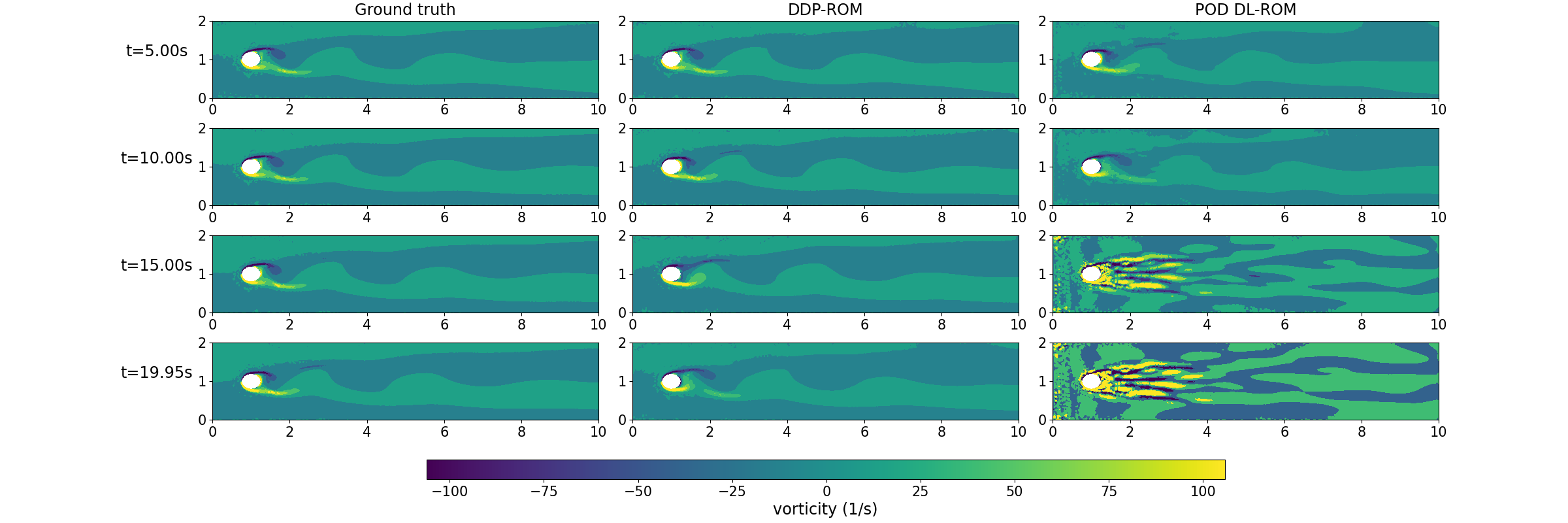}
    \caption{Vorticity snapshots for a trajectory in the test set with $\alpha_{\text{in}}=0.05$ rad and $\gamma_{\text{in}}=-5.99$ m/s. Left column is the ground truth, middle column the DDP-ROM prediction, and right column is DL-ROM prediction.}
    \label{fig:vort_t11}
\end{figure}

In Figure \ref{fig:vmag_t11} and \ref{fig:vort_t11} we show the predictions of the velocity magnitude and vorticity, respectively, of DDP-ROM and DL-ROM compared to the ground truth for a test trajectory $\alpha_{\text{in}}=0.05$ rad and $\gamma_{\text{in}}=5.99$ m/s. Both DDP-ROM and DL-ROM are capable of predicting the velocity and vorticity quite well for $t\leq10$, i.e., the training horizon. However, in the case of time extrapolation ($t>10$), the predictions of DL-ROM are inaccurate and contain many unrealistic artifacts, resulting in a prediction error of $0.47$. Conversely, the predictions of DDP-ROM remain reasonably accurate, although with a little phase shift, resulting in a prediction error of $0.15$.

\begin{figure}[h]
    \centering
    \includegraphics[width=\textwidth]{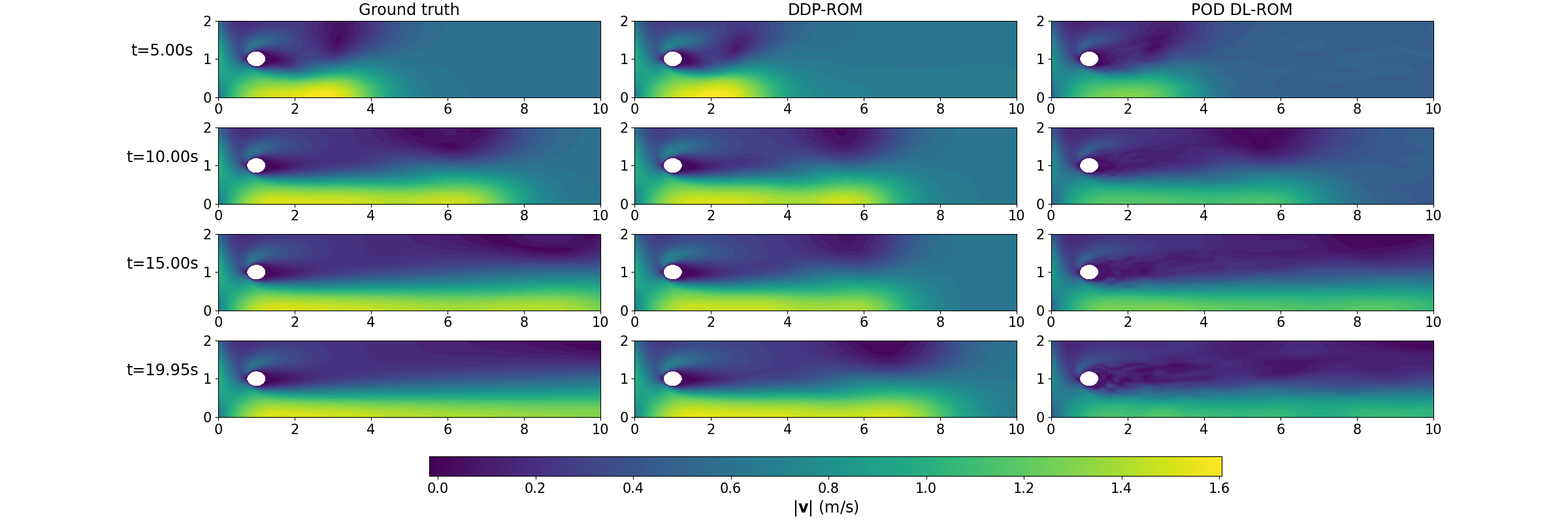}
    \caption{Velocity magnitude snapshots for a trajectory in the test set with $\alpha_{\text{in}}=-0.98$ rad and $\gamma_{\text{in}}=1.05$ m/s. Left column is the ground truth, middle column the DDP-ROM prediction, and right column is DL-ROM prediction.}
    \label{fig:vmag_t28}
\end{figure}

\begin{figure}[h]
    \centering
    \includegraphics[width=\textwidth]{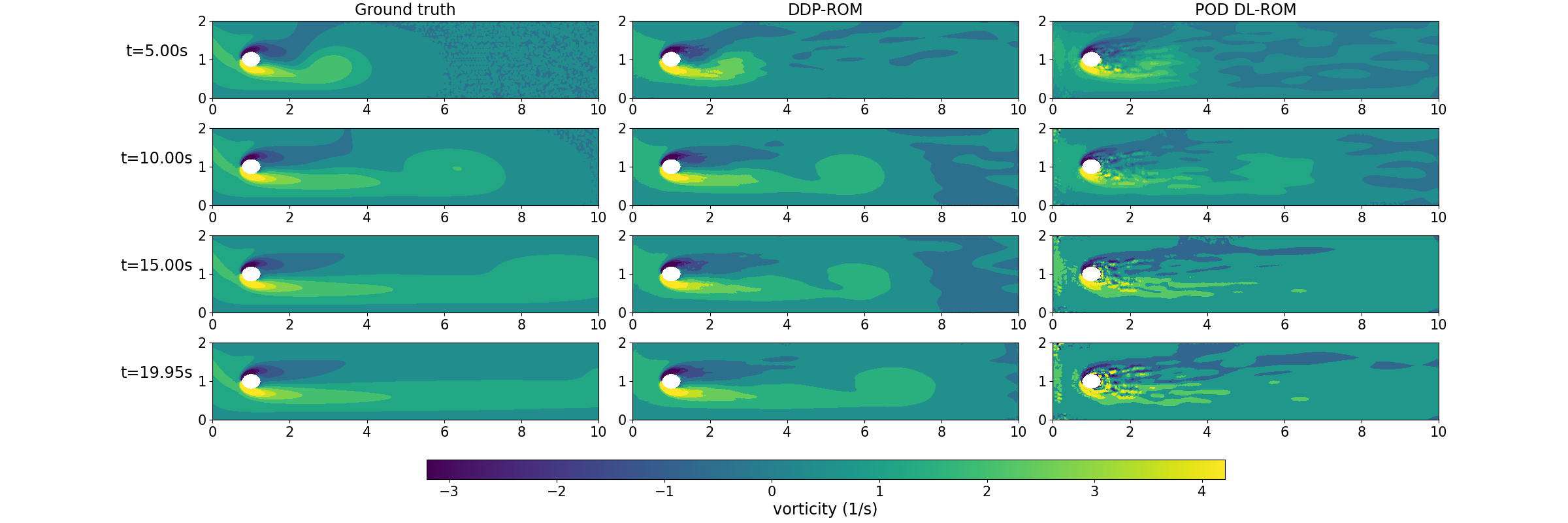}
    \caption{Vorticity snapshots for a trajectory in the test set with $\alpha_{\text{in}}=-0.98$ rad and $\gamma_{\text{in}}=1.05$ m/s. Left column is the ground truth, middle column the DDP-ROM prediction, and right column is DL-ROM prediction.}
    \label{fig:vort_t28}
\end{figure}

In Figure \ref{fig:vmag_t28} and \ref{fig:vort_t28} we show the predictions of the velocity magnitude and vorticity, respectively, of DDP-ROM and DL-ROM compared to the ground truth for a trajectory with $\alpha_{\text{in}}=-0.98$ rad and $\gamma_{\text{in}}=1.05$ m/s. For both DDP-ROM and DL-ROM, this parameter combination yielded the largest mean relative errors across all trajectories in the testing set. DL-ROM underestimates the velocity magnitude and produces much rougher predictions. As a result, the vorticity is not predicted well. On the other side, DDP-ROM is able to predict the velocity features well. However, DDP-ROM struggles at the right-hand side of the domain for $t>10$. This happens because for low velocities for $t<10$ the flow has not reached the right-hand side of the domain yet, making time extrapolation challenging. 

\begin{figure}[h!]
    \centering
    \includegraphics[width=0.8\textwidth]{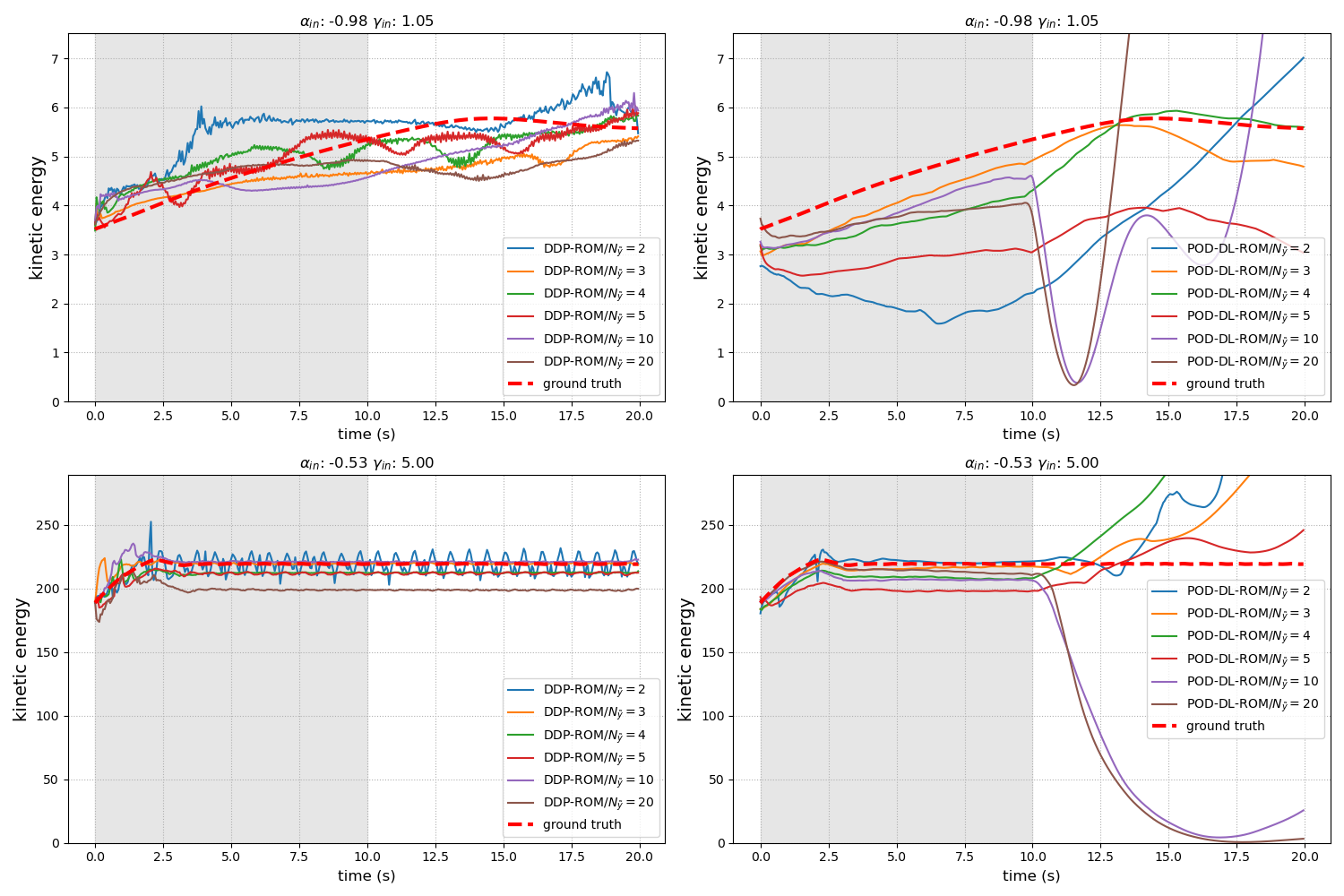}
    \caption{Evolution over time of the kinetic energy of various model predictions against the ground truth for two trajectories in the test set. On top, there is a trajectory with $\alpha_{\text{in}}=-0.98$ rad and $\gamma_{\text{in}}=1.05$ m/s, and at the bottom there is a trajectory with $\alpha_{\text{in}}=-0.53$ rad and $\gamma_{\text{in}}=5.00$ m/s. The gray background highlights the time-interpolation regime for $t\leq10$ s, while the white background the time-extrapolation regime for $t>10$ s.}
    \label{fig:ke_time}
\end{figure}

In Figure \ref{fig:ke_time} we show the evolution over time of the kinetic energy (see \eqref{eq:e_kin}) of DDP-ROM and DL-ROM predictions for different latent dimensions of the inv-AE. and for two different trajectories with parameters $\alpha_{\text{in}}=-0.98$ rad, $\gamma_{\text{in}}=1.05$ m/s, and $\alpha_{\text{in}}=-0.53$ rad and $\gamma_{\text{in}}=5.00$ m/s, respectively. It worth highlighting the different scale of the kinetic energy and in the duration of the transient. In particular, note that for small $\gamma_{\text{in}}$, the training data only contains transients and does not incorporate the steady state, while for larger $\gamma_{\text{in}}$ after an initial peak the kinetic energy becomes constant. When comparing the DDP-ROM and DL-ROM predictions, we observe that the quality of the DL-ROM predictions rapidly deteriorate when extrapolating in time, especially for larger latent dimensions of the inv-AE. The DDP-ROM predictions, on the other hand, remain consistent, even when extrapolating in time, and show stable kinetic energy after the transient has passed.

\begin{figure}[h!]
    \centering
    \includegraphics[width=0.8\textwidth]{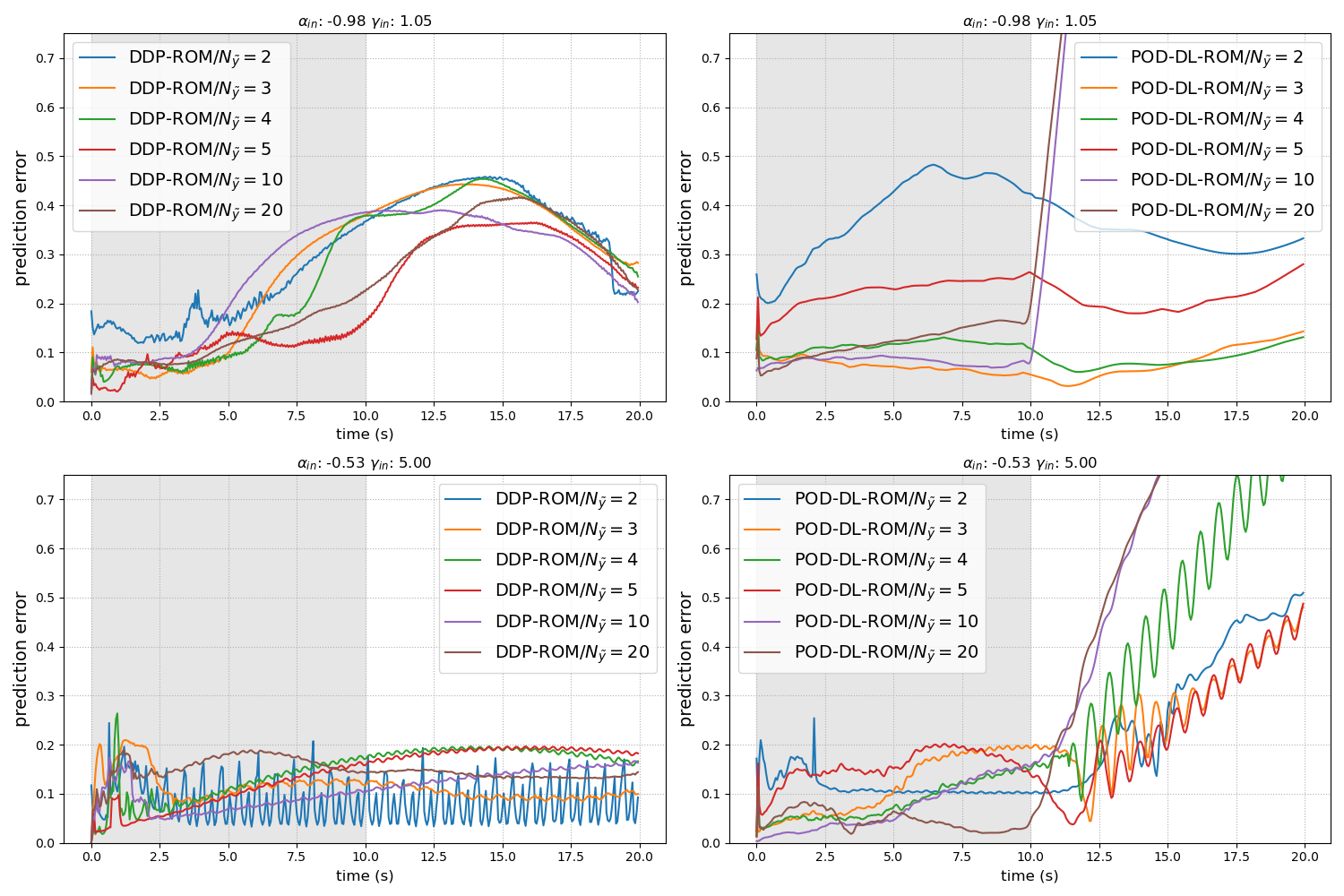}
    \caption{Prediction error over time of various model predictions with an inv-AE architecture for two trajectories in the test set. On top, there is a trajectory with $\alpha_{\text{in}}=-0.98$ rad and $\gamma_{\text{in}}=1.05$ m/s, and at the bottom there is a trajectory with $\alpha_{\text{in}}=1.00$ rad and $\gamma_{\text{in}}=4.68$ m/s. The gray background highlights the time-interpolation regime for $t\leq10$ s, while the white background the time-extrapolation regime for $t>10$ s.}
    \label{fig:err_time}
\end{figure}

In Figure \ref{fig:err_time} we show the relative error over time for the same trajectories as in Figure \ref{fig:ke_time}. For the DDP-ROM, the peaks in relative error approximately match the peaks in kinetic energy in Figure \ref{fig:ke_time}. Next, as the system approaches a steady state, the error decreases again, which is made possible by the spatio-temporal generation of DDP-ROM rather than through autoregression forecasting. The error of DL-ROM predictions, instead, tends to increase in time for $t>10$, where the effect is again stronger in models with a larger latent space.

\begin{figure}[h!]
    \centering
    \includegraphics[width=0.8\textwidth]{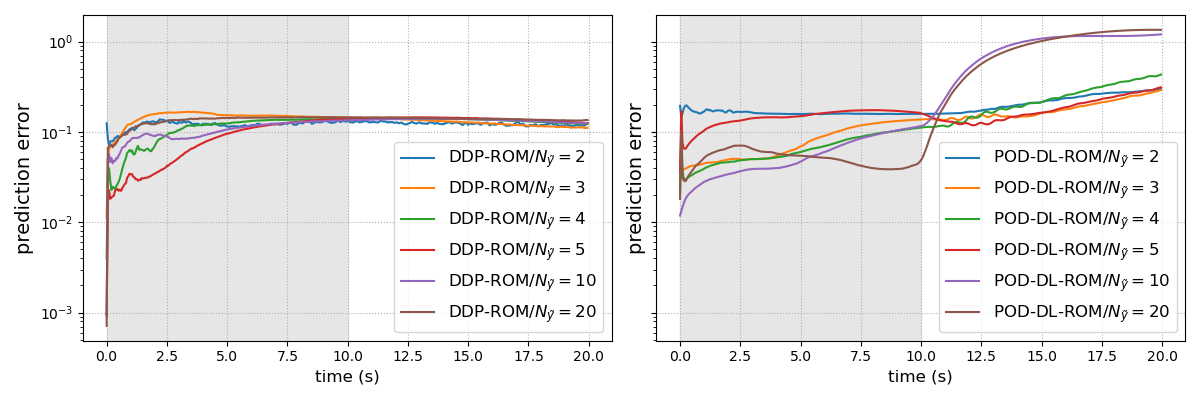}
    \caption{Mean relative error over time of various model predictions for two trajectories averaged over the entire test set. The gray background highlights the time-interpolation regime for $t\leq10$ s, while the white background the time-extrapolation regime for $t>10$ s.}
    \label{fig:err_mean_time}
\end{figure}

Eventually, in Figure \ref{fig:err_mean_time} we depict the average relative error over all trajectories int the training dataset for different latent dimensions of the inv-AE. DDP-ROM is able to keep the relative error stable in time after the initial transient period. Conversely, while DL-ROMs may outperform some of the DDPMs for $t\leq 10$ s for a sufficiently large latent space, the error rapidly increases during time extrapolation ($t>10$ s).

\subsubsection{Uncertainty quantification with DDP-ROM}
\begin{figure}[h!]
    \centering
    \includegraphics[width=\textwidth]{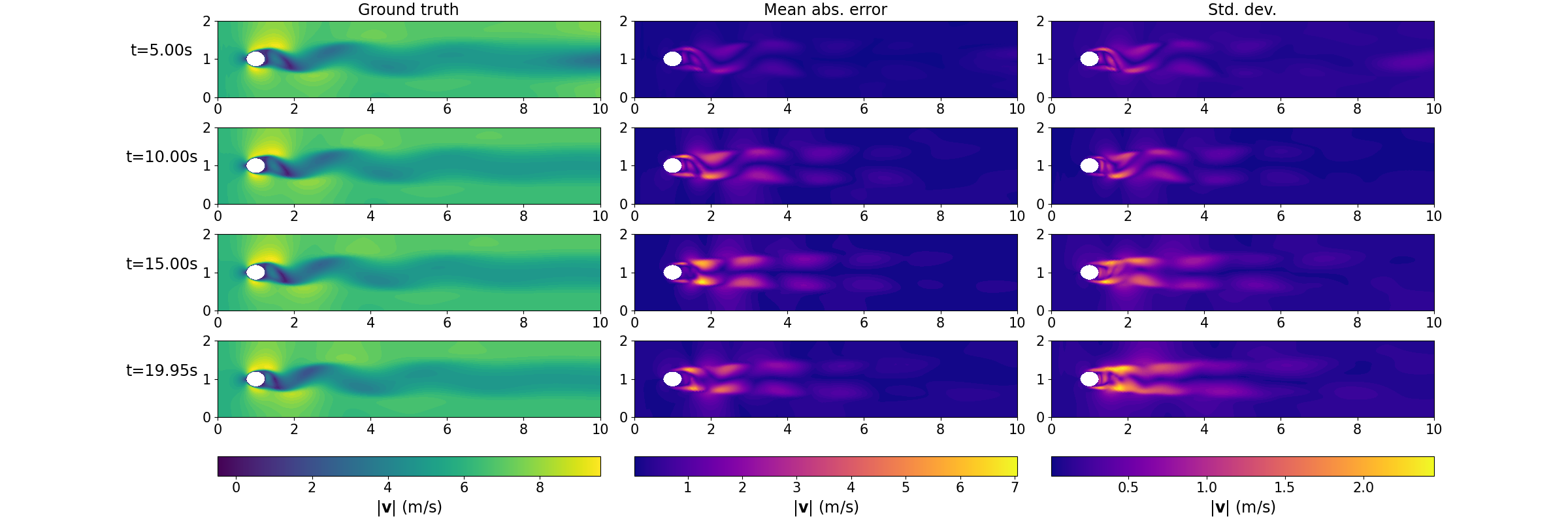}
    \caption{Velocity magnitude snapshot predictions produced by a DDP-ROM with inv-AE architecture with a latent dimension of 5. Parameters are $\alpha_{\text{in}}=0.05$ rad and $\gamma_{\text{in}}=-5.99$ m/s. Left column is the ground truth, middle the mean absolute error and the right column the standard deviation, computed over 10 DDP-ROM predictions.}
    \label{fig:vmag_t11_var}
\end{figure}

\begin{figure}[h!]
    \centering
    \includegraphics[width=\textwidth]{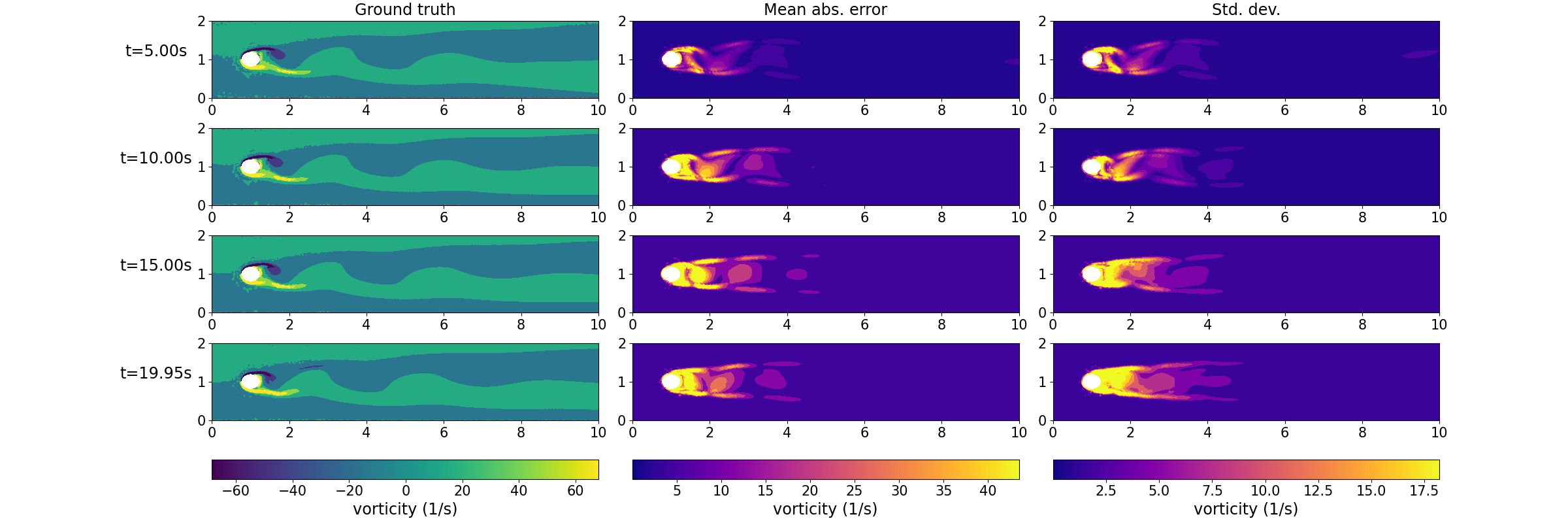}
    \caption{Vorticity snapshot predictions produced by a DDP-ROM with inv-AE architecture with a latent dimension of 5. Parameters are $\alpha_{\text{in}}=0.05$ rad and $\gamma_{\text{in}}=-5.99$ m/s. Left column is the ground truth, middle the mean absolute error and the right column the standard deviation, computed over 10 DDP-ROM predictions.}
    \label{fig:vort_t11_var}
\end{figure}

\begin{figure}[h!]
    \centering
    \includegraphics[width=\textwidth]{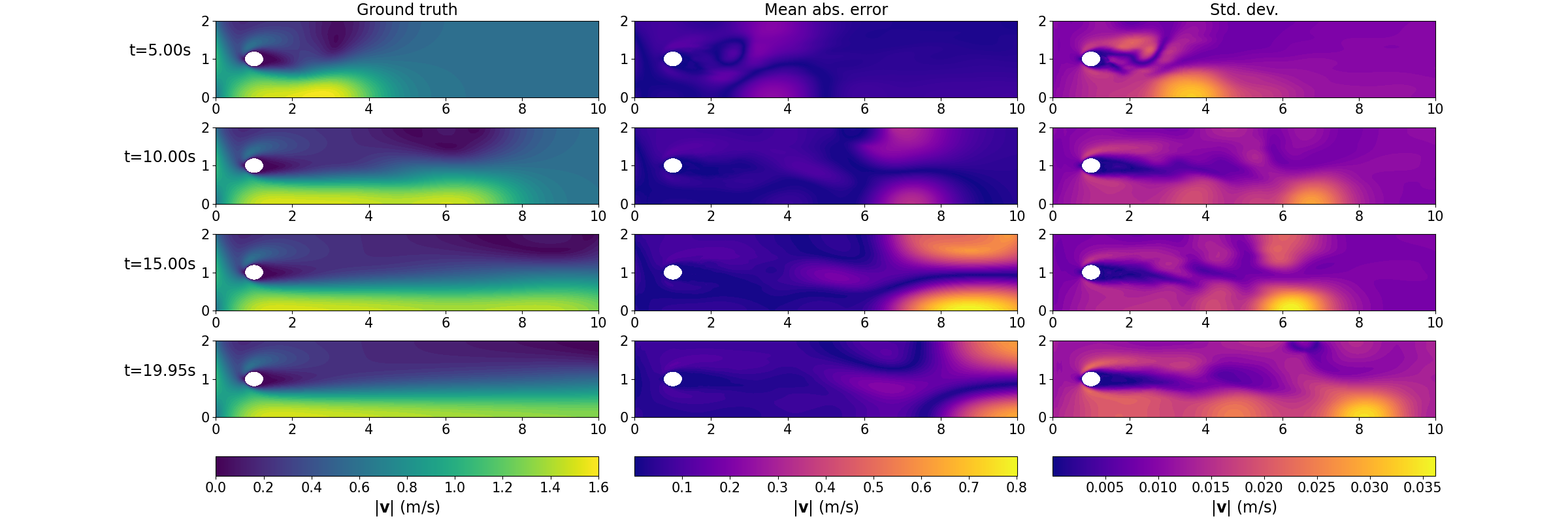}
    \caption{Velocity magnitude snapshot predictions produced by a DDP-ROM with inv-AE architecture with a latent dimension of 5. Parameters are $\alpha_{\text{in}}=-0.98$ rad and $\gamma_{\text{in}}=1.05$ m/s. Left column is the ground truth, middle the mean absolute error and the right column the standard deviation, computed over 10 DDP-ROM predictions.}
    \label{fig:vmag_t28_var}
\end{figure}

\begin{figure}[h!]
    \centering
    \includegraphics[width=\textwidth]{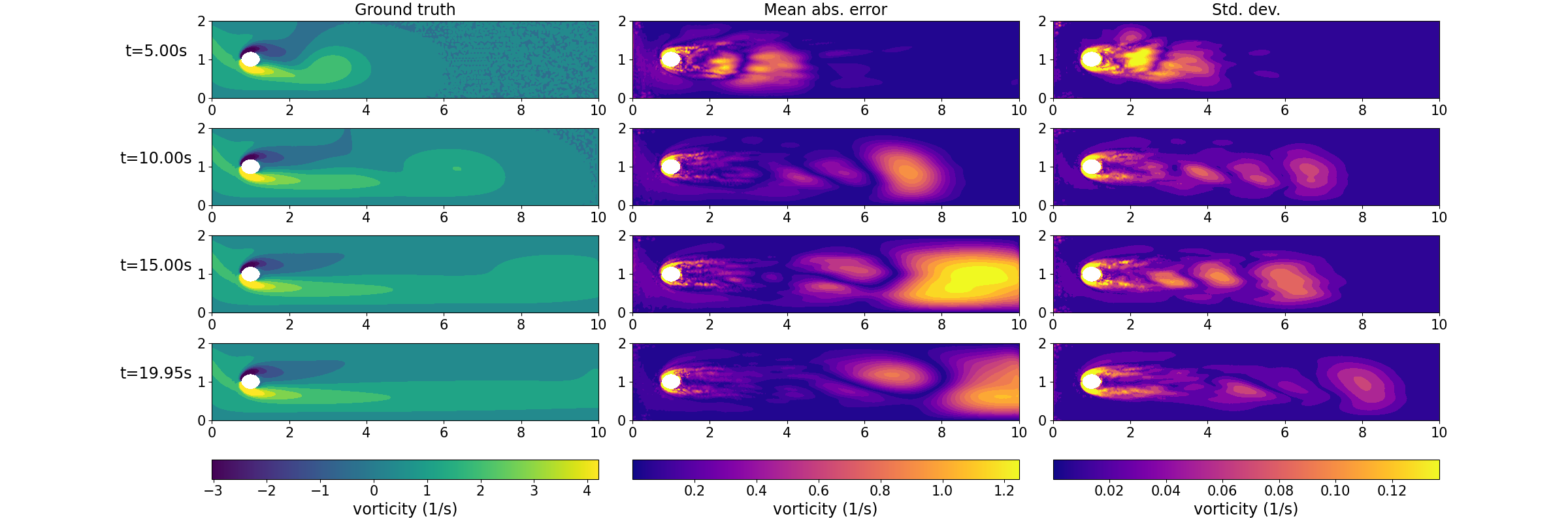}
    \caption{Vorticity snapshot predictions produced by a DDP-ROM with inv-AE architecture with a latent dimension of 5. Parameters are $\alpha_{\text{in}}=-0.98$ rad and $\gamma_{\text{in}}=1.05$ m/s. Left column is the ground truth, the middle column is the mean absolute error and the right column the standard deviation, computed over 10 DDP-ROM predictions.}
    \label{fig:vort_t28_var}
\end{figure}

DDPMs are probabilistic models and, as such, for the same conditioning input $\boldsymbol{\rho}$ and initial condition provided, they can produce different predictions. From multiple of these predictions we can quantify uncertainties.  In Figure \ref{fig:vmag_t11_var} and \ref{fig:vort_t11_var}  we show the mean absolute error and the standard deviation for the velocity magnitude and vorticity, respectively, for $\alpha_{\text{in}}=0.05$ rad and $\gamma_{\text{in}}=-5.99$ m/s, while n Figure \ref{fig:vmag_t28_var} and \ref{fig:vort_t28_var}  we show the mean absolute error and the standard deviation for $\alpha_{\text{in}}=-0.98$ rad and small $\gamma_{\text{in}}=1.05$ m/s. Mean and standard deviation are obtained from 10 DDP-ROM predictions.
Although there is not an exact match between the standard deviation of the predictions and the mean absolute error, the spatial locations where the standard deviation is larger correspond with  regions of higher mean absolute error. In the first case, due to the higher flow velocity, the flow quickly reaches the steady state, making the forecasting problem simpler, while in the second case the transient dominates the flow evolution. In addition, as discussed previously, the time $t>10$ is not captured in the training horizon, which means that for low flow velocities, there is no data on the entire transient, resulting in a larger extrapolation error. As expected, the uncertainties in the first case are smaller than in the second case and DDP-ROM is able to capture this difference.

\begin{figure}[h!]
    \centering
    \includegraphics[width=0.8\textwidth]{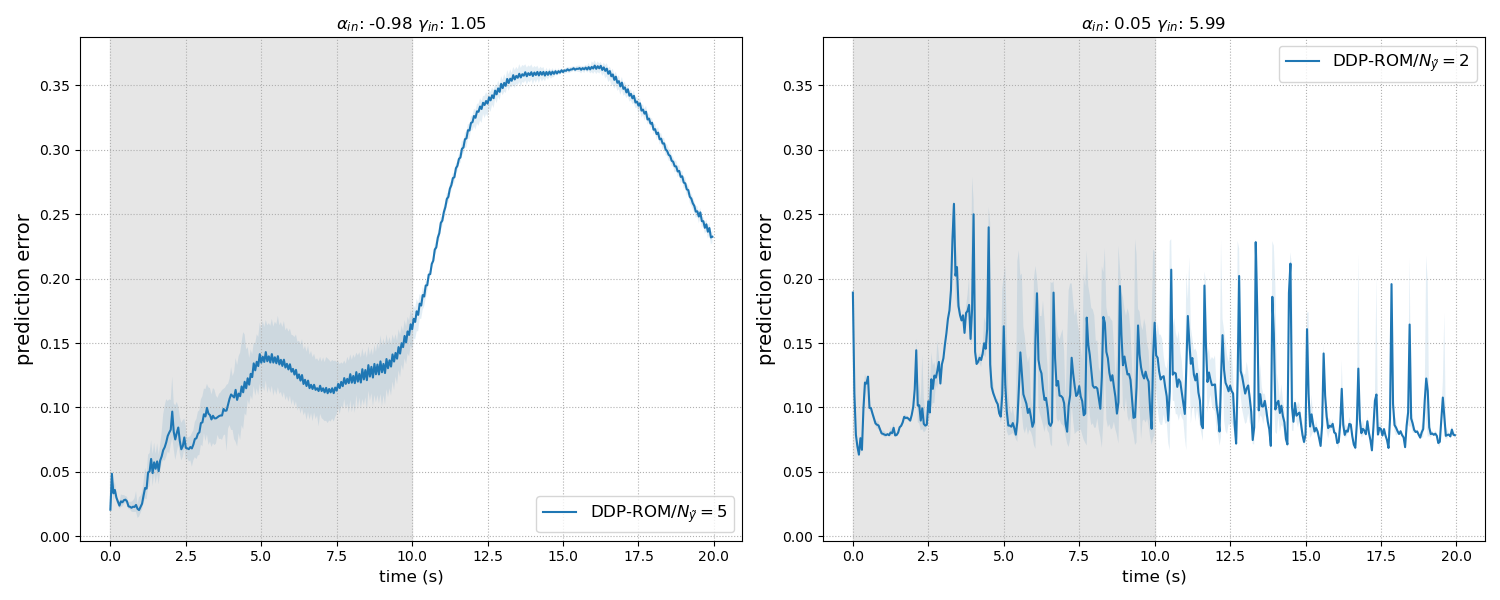}
    \caption{Relative error over time for DDP-ROM predictions with inv-AE architecture with a latent dimension of 2 for two trajectories. The light blue shaded region indicates one standard deviation.}
    \label{fig:err_t_z2}
\end{figure}

\begin{figure}[h!]
    \centering
    \includegraphics[width=0.8\textwidth]{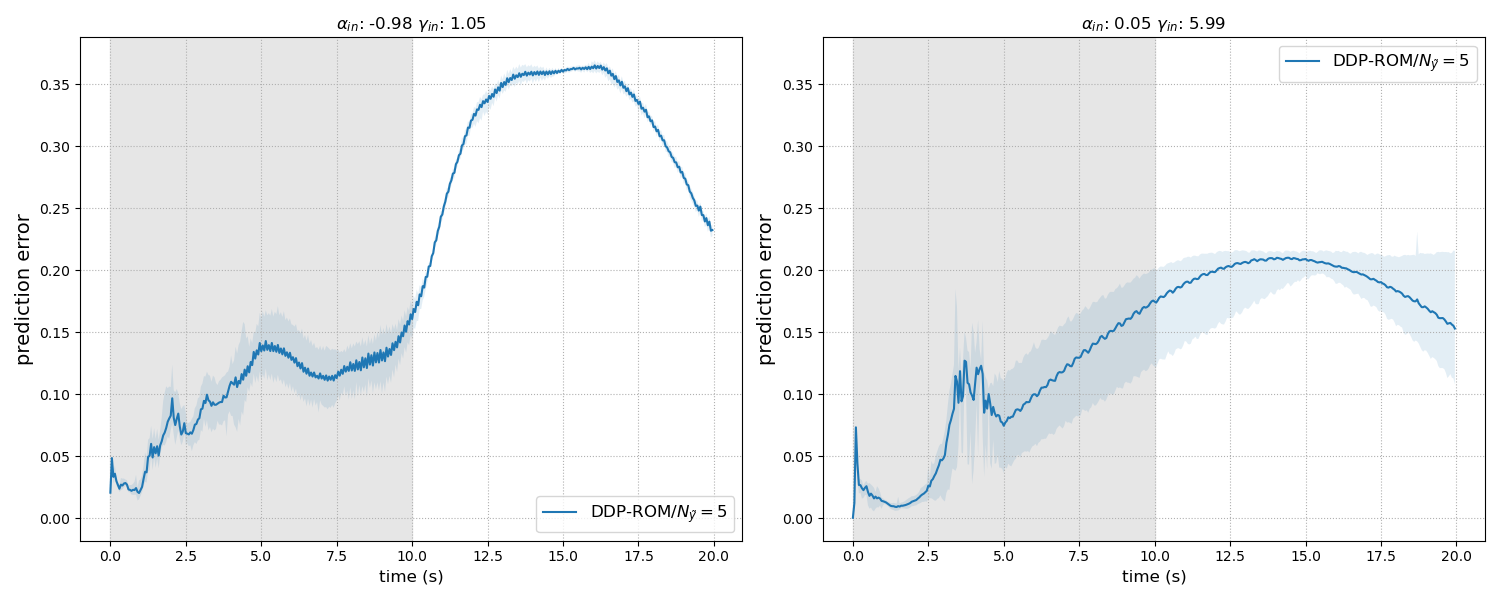}
    \caption{Relative error over time for DDP-ROM predictions with inv-AE architecture with a latent dimension of 5 for two trajectories. The light blue shaded region indicates one standard deviation.}
    \label{fig:err_t_z5}
\end{figure}

\begin{figure}[h!]
    \centering
    \includegraphics[width=0.8\textwidth]{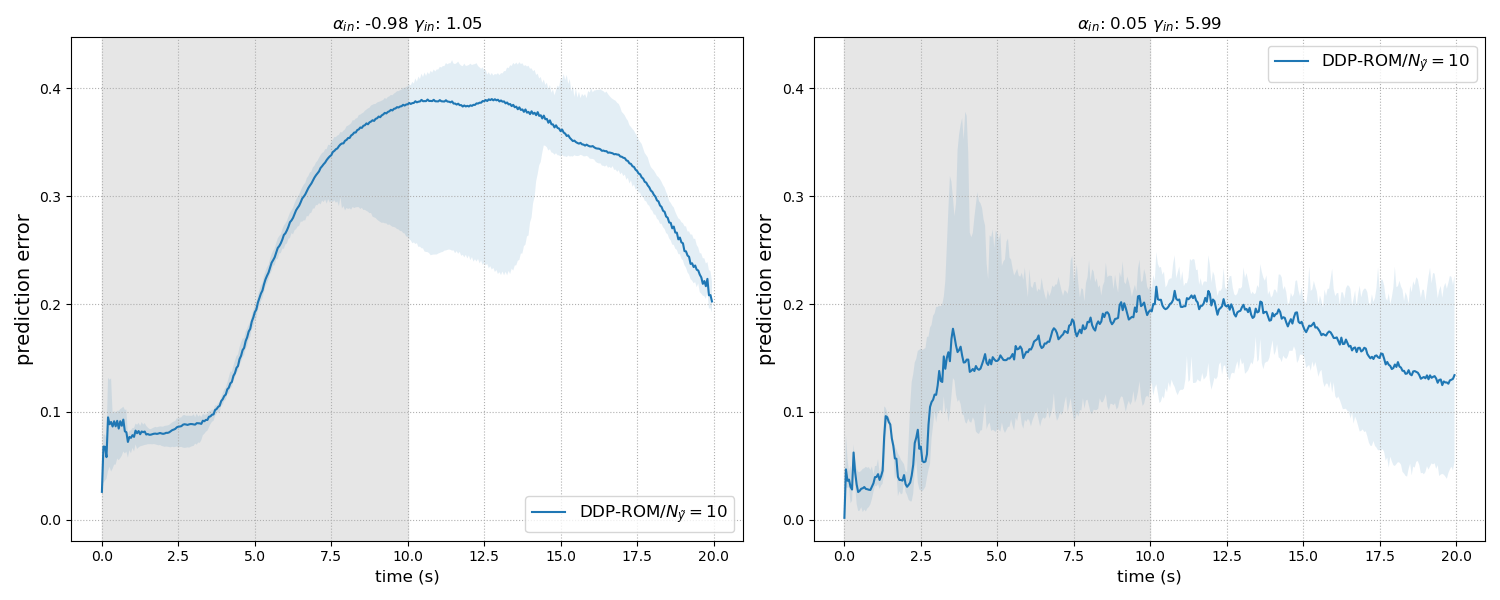}
    \caption{Relative error over time for DDP-ROM predictions with inv-AE architecture with a latent dimension of 10 for two trajectories. The light blue shaded region indicates one standard deviation.}
    \label{fig:err_t_z10}
\end{figure}

Eventually, it is worth highlighting that the variability of the predictions also depends on the dimension of the latent space. Figures \ref{fig:err_t_z2}, \ref{fig:err_t_z5}, and \ref{fig:err_t_z10} show the variation in the relative error over time for a DDP-ROM with inv-AE architecture at various latent dimensions. In particular, a larger latent space dimension leads to a larger spread in the error.
\clearpage

\section{Conclusion}\label{sec:conclusion}
In this work, we proposed a DDPM-based generative ROM, named DDP-ROM, that learns surrogate models of parameterized and high-dimensional dynamical systems. DDP-ROM \emph{(i)} learns a reduced space using POD and an AE in which it embeds a DDPM to mitigate the cost of generating high-dimensional solutions, \emph{(ii)} can make parametric forecasts by using classifier-free condition generation from the latent DDPM, and \emph{(iii)} generates spatio-temporal (latent) trajectories to improve temporal consistency of the solution and enable time extrapolation.

Our findings demonstrate that DDP-ROM can generate coherent and accurate solutions of a parameterized fluid flow around an obstacle with an inference speed that is 400 times faster than the FOM. DDP-ROM produced temporally consistent results across a range of parameters and in time-interpolation and time-extrapolation tasks. Notably, the DDP-ROM can effectively handle unseen parameter values of $\boldsymbol{\rho}$, showing that DDP-ROM is an effective approach for parameterized dynamical systems.  In addition, the probabilistic nature of DDP-ROM can be used to estimate uncertainties of the generated solutions, opening the door towards building uncertainty-aware ROMs.

\section*{Acknowledgments}
The authors would like to express their gratitude to Agata Sowa for their work on an earlier version of this project. \\
NB acknowledges the Project “Reduced Order Modeling and Deep Learning for the real-time approximation of PDEs (DREAM)” (Starting Grant No. FIS00003154), funded by the Italian Science Fund (FIS) - Ministero dell'Università e della Ricerca. NB is a member of the National Group of Scientific Computing (GNCS). 

\bibliographystyle{unsrt}  
\bibliography{references}  

@article{Ma2025Efficient,
  author={Ma, Zhiyuan and Zhang, Yuzhu and Jia, Guoli and Zhao, Liangliang and Ma, Yichao and Ma, Mingjie and Liu, Gaofeng and Zhang, Kaiyan and Ding, Ning and Li, Jianjun and Zhou, Bowen},
  journal={IEEE Transactions on Pattern Analysis and Machine Intelligence}, 
  title={Efficient Diffusion Models: A Comprehensive Survey From Principles to Practices}, 
  year={2025},
  volume={47},
  number={9},
  pages={7506-7525},
  doi={10.1109/TPAMI.2025.3569700}}

@inproceedings{Ohlberger2016Reduced,
    author = {Ohlberger, M and Rave, S.},
    title = {Reduced Basis Methods: Success, Limitations and Future Challenges},
    booktitle = {Proceedings Of The Conference Algoritmy},
    pages = {1--12},
    year = {2016}
}

@article{Sirovich1987Turbulence,
 ISSN = {0033569X, 15524485},
 author = {Lawrence Sirovich},
 journal = {Quarterly of Applied Mathematics},
 number = {3},
 pages = {561--571},
 publisher = {Brown University},
 title = {Turbulence and the dynamics of coherent structures. Part I: Coherent structures},
 volume = {45},
 year = {1987}
}

@article{li2017extended,
    author = {Li, Qianxiao and Dietrich, Felix and Bollt, Erik M. and Kevrekidis, Ioannis G.},
    title = {Extended dynamic mode decomposition with dictionary learning: A data-driven adaptive spectral decomposition of the {Koopman} operator},
    journal = {Chaos: An Interdisciplinary Journal of Nonlinear Science},
    volume = {27},
    number = {10},
    pages = {103111},
    year = {2017},
    month = {10},
    issn = {1054-1500},
    doi = {10.1063/1.4993854},
}

@article{SHARMA2024116865,
title = {Lagrangian operator inference enhanced with structure-preserving machine learning for nonintrusive model reduction of mechanical systems},
journal = {Computer Methods in Applied Mechanics and Engineering},
volume = {423},
pages = {116865},
year = {2024},
issn = {0045-7825},
doi = {https://doi.org/10.1016/j.cma.2024.116865},
author = {Harsh Sharma and David A. Najera-Flores and Michael D. Todd and Boris Kramer},
}

@article{Aubry_Holmes_Lumley_Stone_1988,
    title={The dynamics of coherent structures in the wall region of a turbulent boundary layer},
    volume={192}, 
    DOI={10.1017/S0022112088001818}, 
    journal={Journal of Fluid Mechanics}, 
    author={Aubry, Nadine and Holmes, Philip and Lumley, John L. and Stone, Emily}, 
    year={1988}, 
    pages={115–173}
}

@article{scikit-learn,
  title={Scikit-learn: Machine Learning in {P}ython},
  author={Pedregosa, F. and Varoquaux, G. and Gramfort, A. and Michel, V.
          and Thirion, B. and Grisel, O. and Blondel, M. and Prettenhofer, P.
          and Weiss, R. and Dubourg, V. and Vanderplas, J. and Passos, A. and
          Cournapeau, D. and Brucher, M. and Perrot, M. and Duchesnay, E.},
  journal={Journal of Machine Learning Research},
  volume={12},
  pages={2825--2830},
  year={2011}
}

@article{GILL2026114718,
title = {Fast prediction of plasma instabilities with sparse-grid-accelerated optimized dynamic mode decomposition},
journal = {Journal of Computational Physics},
volume = {553},
pages = {114718},
year = {2026},
issn = {0021-9991},
doi = {https://doi.org/10.1016/j.jcp.2026.114718},
url = {https://www.sciencedirect.com/science/article/pii/S0021999126000689},
author = {Kevin Gill and Ionuţ-Gabriel Farcaş and Silke Glas and Benjamin J. Faber},
}

@article{Hawkins2012,
author = {Hawkins-Daarud, Andrea and van der Zee, Kristoffer G. and Tinsley Oden, J.},
title = {Numerical simulation of a thermodynamically consistent four-species tumor growth model},
journal = {International Journal for Numerical Methods in Biomedical Engineering},
volume = {28},
number = {1},
pages = {3-24},
doi = {https://doi.org/10.1002/cnm.1467},
url = {https://onlinelibrary.wiley.com/doi/abs/10.1002/cnm.1467},
year = {2012}
}

@article{MARCICKI2013310,
title = {Design and parametrization analysis of a reduced-order electrochemical model of graphite/{LiFePO4} cells for {SOC}/{SOH} estimation},
journal = {Journal of Power Sources},
volume = {237},
pages = {310-324},
year = {2013},
issn = {0378-7753},
doi = {https://doi.org/10.1016/j.jpowsour.2012.12.120},
url = {https://www.sciencedirect.com/science/article/pii/S0378775313000694},
author = {James Marcicki and Marcello Canova and A. Terrence Conlisk and Giorgio Rizzoni},
}

@article{Mengaldo2019,
    author = {Mengaldo, Gianmarco and Wyszogrodzki, Andrzej and Diamantakis, Michail and Lock, Sarah-Jane and Giraldo, Francis X. and Wedi, Nils P.},
    title = {Current and Emerging Time-Integration Strategies in Global Numerical Weather and Climate Prediction},
    journal = {Archives of Computational Methods in Engineering},
    year = {2019},
    volume = {26},
    pages = {663–684},
    doi = {https://doi.org/10.1007/s11831-018-9261-8}
}

@article{ZHANG2022120081,
title = {Fully analytical model of heating networks for integrated energy systems},
journal = {Applied Energy},
volume = {327},
pages = {120081},
year = {2022},
issn = {0306-2619},
doi = {https://doi.org/10.1016/j.apenergy.2022.120081},
url = {https://www.sciencedirect.com/science/article/pii/S0306261922013381},
author = {Suhan Zhang and Wei Gu and Xiao-ping Zhang and Hai Lu and Shuai Lu and Ruizhi Yu and Haifeng Qiu}}

@ARTICLE{Paranjape2013,
  author={Paranjape, Aditya A. and Guan, Jinyu and Chung, Soon-Jo and Krstic, Miroslav},
  journal={IEEE Transactions on Robotics}, 
  title={PDE Boundary Control for Flexible Articulated Wings on a Robotic Aircraft}, 
  year={2013},
  volume={29},
  number={3},
  pages={625-640},
  doi={10.1109/TRO.2013.2240711}}

@article{botteghi2022deep,
  title={Deep kernel learning of dynamical models from high-dimensional noisy data},
  author={Botteghi, Nicol{\`o} and Guo, Mengwu and Brune, Christoph},
  journal={Scientific Reports},
  volume={12},
  number={1},
  pages={21530},
  year={2022},
  publisher={Nature Publishing Group UK London}
}

@book{sudret2000stochastic,
  title={Stochastic finite element methods and reliability: a state-of-the-art report},
  author={Sudret, Bruno and Der Kiureghian, Armen},
  year={2000},
  publisher={Department of Civil and Environmental Engineering, University of California}
}

@book{quarteroni2015reduced,
  title={Reduced basis methods for partial differential equations: an introduction},
  author={Quarteroni, Alfio and Manzoni, Andrea and Negri, Federico},
  volume={92},
  year={2015},
  publisher={Springer}
}

@book{brunton2019data,
  title={Data-driven science and engineering: Machine learning, dynamical systems, and control},
  author={Brunton, Steven L and Kutz, J Nathan},
  year={2019},
  publisher={Cambridge University Press}
}

@article{schmid2010dynamic,
  title={Dynamic mode decomposition of numerical and experimental data},
  author={Schmid, Peter J},
  journal={Journal of Fluid Mechanics},
  volume={656},
  pages={5--28},
  year={2010},
  publisher={Cambridge University Press}
}

@article{guo2022bayesian,
  title={Bayesian operator inference for data-driven reduced-order modeling},
  author={Guo, Mengwu and McQuarrie, Shane A and Willcox, Karen E},
  journal={Computer Methods in Applied Mechanics and Engineering},
  volume={402},
  pages={115336},
  year={2022},
  publisher={Elsevier}
}

@article{peherstorfer2016data,
  title={Data-driven operator inference for nonintrusive projection-based model reduction},
  author={Peherstorfer, Benjamin and Willcox, Karen},
  journal={Computer Methods in Applied Mechanics and Engineering},
  volume={306},
  pages={196--215},
  year={2016},
  publisher={Elsevier}
}

@article{qian2020lift,
  title={{Lift \& learn}: Physics-informed machine learning for large-scale nonlinear dynamical systems},
  author={Qian, Elizabeth and Kramer, Boris and Peherstorfer, Benjamin and Willcox, Karen},
  journal={Physica D: Nonlinear Phenomena},
  volume={406},
  pages={132401},
  year={2020},
  publisher={Elsevier}
}

@article{bakarji2022discovering,
    author = {Bakarji, Joseph and Champion, Kathleen and Nathan Kutz, J. and Brunton, Steven L.},
    title = {Discovering governing equations from partial measurements with deep delay autoencoders},
    journal = {Proceedings of the Royal Society A: Mathematical, Physical and Engineering Sciences},
    volume = {479},
    number = {2276},
    pages = {20230422},
    year = {2023},
    month = {08},
    issn = {1364-5021},
    doi = {10.1098/rspa.2023.0422},
}

@article{brunton2016discovering,
  title={Discovering governing equations from data by sparse identification of nonlinear dynamical systems},
  author={Brunton, Steven L and Proctor, Joshua L and Kutz, J Nathan},
  journal={Proceedings of the National Academy of Sciences},
  volume={113},
  number={15},
  pages={3932--3937},
  year={2016},
  publisher={National Acad Sciences}
}

@article{champion2019data,
  title={Data-driven discovery of coordinates and governing equations},
  author={Champion, Kathleen and Lusch, Bethany and Kutz, J Nathan and Brunton, Steven L},
  journal={Proceedings of the National Academy of Sciences},
  volume={116},
  number={45},
  pages={22445--22451},
  year={2019},
  publisher={National Acad Sciences}
}

@article{fresca2021comprehensive,
  title={A comprehensive deep learning-based approach to reduced order modeling of nonlinear time-dependent parametrized {PDE}s},
  author={Fresca, Stefania and Ded\'{e}, Luca and Manzoni, Andrea},
  journal={Journal of Scientific Computing},
  volume={87},
  pages={1--36},
  year={2021},
  publisher={Springer}
}

@article{fresca2022pod,
  title={{POD-DL-ROM}: Enhancing deep learning-based reduced order models for nonlinear parametrized PDEs by proper orthogonal decomposition},
  author={Fresca, Stefania and Manzoni, Andrea},
  journal={Computer Methods in Applied Mechanics and Engineering},
  volume={388},
  pages={114181},
  year={2022},
  publisher={Elsevier}
}

@article{otto2019linearly,
  title={Linearly recurrent autoencoder networks for learning dynamics},
  author={Otto, Samuel E and Rowley, Clarence W},
  journal={SIAM Journal on Applied Dynamical Systems},
  volume={18},
  number={1},
  pages={558--593},
  year={2019},
  publisher={SIAM}
}

@article{guo2019data,
  title={Data-driven reduced order modeling for time-dependent problems},
  author={Guo, Mengwu and Hesthaven, Jan S},
  journal={Computer Methods in Applied Mechanics and Engineering},
  volume={345},
  pages={75--99},
  year={2019},
  publisher={Elsevier}
}

@book{goodfellow2016deep,
  title={Deep learning},
  author={Goodfellow, Ian and Bengio, Yoshua and Courville, Aaron},
  year={2016},
  publisher={MIT press}
}

@article{hinton1993autoencoders,
     author = {Hinton, Geoffrey E and Zemel, Richard},
     journal = {Advances in Neural Information Processing Systems},
     editor = {J. Cowan and G. Tesauro and J. Alspector},
     pages = {3--10},
     publisher = {Morgan-Kaufmann},
     title = {Autoencoders, Minimum Description Length and {Helmholtz} Free Energy},
     volume = {6},
     year = {1993}
}

@article{willcox2024role,
  title={The role of computational science in digital twins},
  author={Willcox, Karen and Segundo, Brittany},
  journal={Nature Computational Science},
  volume={4},
  number={3},
  pages={147--149},
  year={2024},
  publisher={Nature Publishing Group}
}

@article{galbally2010non,
  title={Non-linear model reduction for uncertainty quantification in large-scale inverse problems},
  author={Galbally, David and Fidkowski, Krzysztof and Willcox, Karen and Ghattas, Omar},
  journal={International Journal for Numerical Methods in Engineering},
  volume={81},
  number={12},
  pages={1581--1608},
  year={2010},
  publisher={Wiley Online Library}
}

@article{ravindran2000reduced,
  title={A reduced-order approach for optimal control of fluids using proper orthogonal decomposition},
  author={Ravindran, Sivaguru S},
  journal={International Journal for Numerical Methods in Fluids},
  volume={34},
  number={5},
  pages={425--448},
  year={2000},
  publisher={Wiley Online Library}
}

@book{troltzsch2024optimal,
  title={Optimal control of partial differential equations: theory, methods and applications},
  author={Tr{\"o}ltzsch, Fredi},
  volume={112},
  year={2010},
  publisher={American Mathematical Society}
}

@inproceedings{rombach2022high,
  title={High-resolution image synthesis with latent diffusion models},
  author={Rombach, Robin and Blattmann, Andreas and Lorenz, Dominik and Esser, Patrick and Ommer, Bj{\"o}rn},
  booktitle={Proceedings of the IEEE/CVF conference on computer vision and pattern recognition},
  pages={10684--10695},
  year={2022}
}

@article{ajay2022conditional,
  title={Is conditional generative modeling all you need for decision-making?},
  author={Ajay, Anurag and Du, Yilun and Gupta, Abhi and Tenenbaum, Joshua and Jaakkola, Tommi and Agrawal, Pulkit},
  journal={arXiv preprint arXiv:2211.15657},
  year={2022}
}

@article{ho2022classifier,
  title={Classifier-free diffusion guidance},
  author={Ho, Jonathan and Salimans, Tim},
  journal={arXiv preprint arXiv:2207.12598},
  year={2022}
}

@inproceedings{Kingma2014Auto-encodingBayes,
    title = {{Auto-encoding variational bayes}},
    year = {2014},
    booktitle = {2nd International Conference on Learning Representations, ICLR 2014 - Conference Track Proceedings},
    author = {Kingma, Diederik P and Welling, Max},
    arxivId = {1312.6114}
}

@inproceedings{Sohl-Dickstein2015DeepThermodynamics,
    title = 	 {Deep Unsupervised Learning using Nonequilibrium Thermodynamics},
    author = 	 {Sohl-Dickstein, Jascha and Weiss, Eric and Maheswaranathan, Niru and Ganguli, Surya},
    booktitle = 	 {Proceedings of the 32nd International Conference on Machine Learning},
    pages = 	 {2256--2265},
    year = 	 {2015},
    editor = 	 {Bach, Francis and Blei, David},
    volume = 	 {37},
    series = 	 {Proceedings of Machine Learning Research},
    address = 	 {Lille, France},
    month = 	 {07--09 Jul},
    publisher =    {PMLR},
}

@inproceedings{Ho2020DenoisingModels,
     author = {Ho, Jonathan and Jain, Ajay and Abbeel, Pieter},
     booktitle = {Advances in Neural Information Processing Systems},
     editor = {H. Larochelle and M. Ranzato and R. Hadsell and M.F. Balcan and H. Lin},
     pages = {6840--6851},
     publisher = {Curran Associates, Inc.},
     title = {Denoising Diffusion Probabilistic Models},
     volume = {33},
     year = {2020}
}

@inproceedings{Dhariwal2021DiffusionSynthesis,
     author = {Dhariwal, Prafulla and Nichol, Alexander},
     booktitle = {Advances in Neural Information Processing Systems},
     editor = {M. Ranzato and A. Beygelzimer and Y. Dauphin and P.S. Liang and J. Wortman Vaughan},
     pages = {8780--8794},
     publisher = {Curran Associates, Inc.},
     title = {Diffusion Models Beat GANs on Image Synthesis},
     volume = {34},
     year = {2021}
}

@article{Goodfellow2020GenerativeNetworks,
    author = {Goodfellow, Ian and Pouget-Abadie, Jean and Mirza, Mehdi and Xu, Bing and Warde-Farley, David and Ozair, Sherjil and Courville, Aaron and Bengio, Yoshua},
    title = {Generative adversarial networks},
    year = {2020},
    issue_date = {November 2020},
    publisher = {Association for Computing Machinery},
    address = {New York, NY, USA},
    volume = {63},
    number = {11},
    issn = {0001-0782},
    doi = {10.1145/3422622},
    journal = {Commun. ACM},
    month = oct,
    pages = {139–144},
    numpages = {6}
}

@inproceedings{Song2019GenerativeDistribution,
    title = {{Generative modeling by estimating gradients of the data distribution}},
    year = {2019},
    booktitle = {Advances in Neural Information Processing Systems},
    author = {Song, Yang and Ermon, Stefano},
    volume = {32},
    issn = {10495258},
    arxivId = {1907.05600}
}

@InProceedings{Nichol2021ImprovedModels,
  title = 	 {Improved Denoising Diffusion Probabilistic Models},
  author =       {Nichol, Alexander Quinn and Dhariwal, Prafulla},
  booktitle = 	 {Proceedings of the 38th International Conference on Machine Learning},
  pages = 	 {8162--8171},
  year = 	 {2021},
  editor = 	 {Meila, Marina and Zhang, Tong},
  volume = 	 {139},
  series = 	 {Proceedings of Machine Learning Research},
  month = 	 {18--24 Jul},
  publisher =    {PMLR}
}

@inproceedings{Song2020ImprovedModels,
     author = {Song, Yang and Ermon, Stefano},
     booktitle = {Advances in Neural Information Processing Systems},
     editor = {H. Larochelle and M. Ranzato and R. Hadsell and M.F. Balcan and H. Lin},
     pages = {12438--12448},
     publisher = {Curran Associates, Inc.},
     title = {Improved Techniques for Training Score-Based Generative Models},
     volume = {33},
     year = {2020}
}

@InProceedings{Janner2022PlanningSynthesis,
  title = 	 {Planning with Diffusion for Flexible Behavior Synthesis},
  author =       {Janner, Michael and Du, Yilun and Tenenbaum, Joshua and Levine, Sergey},
  booktitle = 	 {Proceedings of the 39th International Conference on Machine Learning},
  pages = 	 {9902--9915},
  year = 	 {2022},
  editor = 	 {Chaudhuri, Kamalika and Jegelka, Stefanie and Song, Le and Szepesvari, Csaba and Niu, Gang and Sabato, Sivan},
  volume = 	 {162},
  series = 	 {Proceedings of Machine Learning Research},
  month = 	 {17--23 Jul},
  publisher =    {PMLR},
}

@article{kadeethum2021framework,
  title={A framework for data-driven solution and parameter estimation of {PDEs} using conditional generative adversarial networks},
  author={Kadeethum, Teeratorn and O’Malley, Daniel and Fuhg, Jan Niklas and Choi, Youngsoo and Lee, Jonghyun and Viswanathan, Hari S and Bouklas, Nikolaos},
  journal={Nature Computational Science},
  volume={1},
  number={12},
  pages={819--829},
  year={2021},
  publisher={Nature Publishing Group US New York}
}

@article{kemna2023reduced,
  title={Reduced order fluid modeling with generative adversarial networks},
  author={Kemna, Mirko and Heinlein, Alexander and Vuik, Cornelis},
  journal={Proceedings in Applied Mathematics and Mechanics},
  volume={23},
  number={1},
  doi = {https://doi.org/10.1002/pamm.202200241},
  pages={e202200241},
  year={2023},
  publisher={Wiley Online Library}
}

@article{coscia2024generative,
  title={Generative adversarial reduced order modelling},
  author={Coscia, Dario and Demo, Nicola and Rozza, Gianluigi},
  journal={Scientific Reports},
  volume={14},
  number={1},
  pages={3826},
  year={2024},
  publisher={Nature Publishing Group UK London}
}

@article{shu2023physics,
  title={A physics-informed diffusion model for high-fidelity flow field reconstruction},
  author={Shu, Dule and Li, Zijie and Farimani, Amir Barati},
  journal={Journal of Computational Physics},
  volume={478},
  pages={111972},
  year={2023},
  publisher={Elsevier}
}

@article{shan2024pird,
  author={Shan, Siming
and Wang, Pengkai
and Chen, Song
and Liu, Jiaxu
and Xu, Chao
and Cai, Shengze},
title={{PiRD}: physics-informed residual diffusion for flow field reconstruction},
journal={Acta Mechanica Sinica},
year={2026},
volume={42},
pages={725259},
issn={1614-3116},
doi={10.1007/s10409-025-25259-x},
}

@article{bastek2024physics,
  title={Physics-Informed Diffusion Models},
  author={Bastek, Jan-Hendrik and Sun, WaiChing and Kochmann, Dennis M},
  journal={arXiv preprint arXiv:2403.14404},
  year={2024}
}

@inproceedings{
xu2022geodiff,
title={GeoDiff: A Geometric Diffusion Model for Molecular Conformation Generation},
author={Minkai Xu and Lantao Yu and Yang Song and Chence Shi and Stefano Ermon and Jian Tang},
booktitle={International Conference on Learning Representations},
year={2022},
}

@article{dureth2023conditional,
  title={Conditional diffusion-based microstructure reconstruction},
  author={D{\"u}reth, Christian and Seibert, Paul and R{\"u}cker, Dennis and Handford, Stephanie and K{\"a}stner, Markus and Gude, Maik},
  journal={Materials Today Communications},
  volume={35},
  pages={105608},
  year={2023},
  publisher={Elsevier}
}

@article{bastek2023inverse,
  title={Inverse design of nonlinear mechanical metamaterials via video denoising diffusion models},
  author={Bastek, Jan-Hendrik and Kochmann, Dennis M},
  journal={Nature Machine Intelligence},
  volume={5},
  number={12},
  pages={1466--1475},
  year={2023},
  publisher={Nature Publishing Group UK London}
}

@inproceedings{zhang2024xddpm,
  title={{XDDPM}: EXPLAINABLE DENOISING DIFFUSION PROBABILISTIC MODEL FOR SCIENTIFIC MODELING},
  author={Zhang, Qianru and Yu, Chenglei and Yan, Yudong and Kuang, Xiangyu and Ma, Yi and Cao, Yuansheng and Yiu, Siu Ming and Wu, Tailin},
  booktitle={ICLR 2024 Workshop on AI4DifferentialEquations In Science},
  year={2024},
}

@inproceedings{shysheya2024conditional,
     author = {Shysheya, Aliaksandra and Diaconu, Cristiana and Bergamin, Federico and Perdikaris, Paris and Hern\'{a}ndez-Lobato, Jos\'{e} Miguel and Turner, Richard E. and Mathieu, Emile},
     booktitle = {Advances in Neural Information Processing Systems},
     doi = {10.52202/079017-0732},
     editor = {A. Globerson and L. Mackey and D. Belgrave and A. Fan and U. Paquet and J. Tomczak and C. Zhang},
     pages = {23246--23300},
     publisher = {Curran Associates, Inc.},
     title = {On conditional diffusion models for PDE simulations},
     volume = {37},
     year = {2024}
}

@InProceedings{botteghi2024recurrentdeepkernellearning,
author="Botteghi, Nicol{\`o}
and Motta, Paolo
and Manzoni, Andrea
and Zunino, Paolo
and Guo, Mengwu",
title="Recurrent Deep Kernel Learning of Dynamical Systems",
booktitle="Physics-Based and Data-Driven Modeling for Digital Twins",
year="2026",
publisher="Springer Nature Singapore",
pages="65--86",
}

@inproceedings{li2020fourier,
  title={Fourier Neural Operator for Parametric Partial Differential Equations},
    author={Zongyi Li and Nikola Borislavov Kovachki and Kamyar Azizzadenesheli and Burigede liu and Kaushik Bhattacharya and Andrew Stuart and Anima Anandkumar},
    booktitle={International Conference on Learning Representations},
    year={2021},
}

@article{kovachki2023neural,
  title={Neural operator: Learning maps between function spaces with applications to {PDEs}},
  author={Kovachki, Nikola and Li, Zongyi and Liu, Burigede and Azizzadenesheli, Kamyar and Bhattacharya, Kaushik and Stuart, Andrew and Anandkumar, Anima},
  journal={Journal of Machine Learning Research},
  volume={24},
  number={89},
  pages={1--97},
  year={2023}
}

@article{cao2023lno,
  title={{LNO}: {Laplace} neural operator for solving differential equations},
  author={Cao, Qianying and Goswami, Somdatta and Karniadakis, George Em},
  journal={arXiv preprint arXiv:2303.10528},
  year={2023}
}

@inproceedings{huang2024diffusionpde,
   author = {Huang, Jiahe and Yang, Guandao and Wang, Zichen and Park, Jeong Joon},
 booktitle = {Advances in Neural Information Processing Systems},
 doi = {10.52202/079017-4140},
 editor = {A. Globerson and L. Mackey and D. Belgrave and A. Fan and U. Paquet and J. Tomczak and C. Zhang},
 pages = {130291--130323},
 publisher = {Curran Associates, Inc.},
 title = {DiffusionPDE: Generative {PDE}-Solving under Partial Observation},
 volume = {37},
 year = {2024}
}

@inproceedings{wei2024generative,
  title={Generative {PDE} Control},
  author={Wei, Long and Hu, Peiyan and Feng, Ruiqi and Du, Yixuan and Zhang, Tao and Wang, Rui and Wang, Yue and Ma, Zhi-Ming and Wu, Tailin},
  booktitle={ICLR 2024 Workshop on AI4DifferentialEquations In Science},
  year={2024},
}

@article{botteghi2026deep,
    title={Deep Invertible Autoencoders for Dimensionality Reduction of Dynamical Systems},
    author={Botteghi, Nicolò and Glas, Silke and Brune, Christoph},
    journal={arXiv preprint arXiv:2603.13496},
    year={2026},
}

@Article{KimCWZ22,
  author   = {Youngkyu Kim and Youngsoo Choi and David Widemann and Tarek Zohdi},
  title    = {A fast and accurate physics-informed neural network reduced order model with shallow masked autoencoder},
  journal  = {Journal of Computational Physics},
  year     = {2022},
  volume   = {451},
  pages    = {110841},
  issn     = {0021-9991},
  doi      = {10.1016/j.jcp.2021.110841},
}

@article{alnaes2015fenics,
  title={The {FEniCS} project version 1.5},
  author={Aln{\ae}s, Martin and Blechta, Jan and Hake, Johan and Johansson, August and Kehlet, Benjamin and Logg, Anders and Richardson, Chris and Ring, Johannes and Rognes, Marie E and Wells, Garth N},
  journal={Archive of Numerical Software},
  volume={3},
  number={100},
  year={2015}
}

@article{chakraborty2021role,
  title={The role of surrogate models in the development of digital twins of dynamic systems},
  author={Chakraborty, Souvik and Adhikari, Sondipon and Ganguli, Ranjan},
  journal={Applied Mathematical Modelling},
  volume={90},
  pages={662--681},
  year={2021},
  publisher={Elsevier}
}

@article{yang2023diffusion,
  title={Diffusion models: A comprehensive survey of methods and applications},
  author={Yang, Ling and Zhang, Zhilong and Song, Yang and Hong, Shenda and Xu, Runsheng and Zhao, Yue and Zhang, Wentao and Cui, Bin and Yang, Ming-Hsuan},
  journal={ACM Computing Surveys},
  volume={56},
  number={4},
  pages={1--39},
  year={2023},
  publisher={ACM New York, NY, USA}
}

@article{price2023gencast,
  title={Probabilistic weather forecasting with machine learning},
  author={Price, Ilan and Sanchez-Gonzalez, Alvaro and Alet, Ferran and Andersson, Tom R and El-Kadi, Andrew and Masters, Dominic and Ewalds, Timo and Stott, Jacklynn and Mohamed, Shakir and Battaglia, Peter and Lam, Remi and Willson, Matthew},
  journal={Nature},
  year={2025},
  volume={637},
  number={8044},
  pages={84-90},
  issn={1476-4687},
  doi={10.1038/s41586-024-08252-9},
  url={https://doi.org/10.1038/s41586-024-08252-9}
}

\clearpage
\appendix

\section{Deep Learning-based ROMs}\label{app:dl-rom}
Deep Learning-based ROM (DL-ROM) \cite{fresca2021comprehensive} is a widely adopted and effective framework for learning ROMs directly from data. DL-ROM consists of an encoder 
\begin{equation}
    \enc(\bm{y}(t_n;\boldsymbol{\rho}_p);\boldsymbol{\theta}_{\enc})\, ,
\end{equation}
and a decoder
\begin{equation}
    \dec(\tilde{\bm{y}}(t_n;\boldsymbol{\rho}_p);\boldsymbol{\theta}_{\dec})\, ,
\end{equation}
which jointly learn a low-dimensional representation $\tilde{\bm{y}}(t;\boldsymbol{\rho}_p)$ of the snapshots and its reconstruction $\hat{\bm{y}}(t_n;\boldsymbol{\rho}_p)$. 
To efficiently predict the system's dynamics, DL-ROM introduces a feedforward neural network 
\begin{equation}
    \xirom(\boldsymbol{\rho}_p, t_n;\boldsymbol{\theta}_{\xirom})\, ,
\end{equation}
with learnable parameters $\boldsymbol{\theta}_{\xi}$, which predicts the latent snapshot $\tilde{\bm{y}}(t_n;\boldsymbol{\rho}_p)$ directly from time $t_n$ and parameter vector $\boldsymbol{\rho}_p$. Approximations of the FOM snapshots are then obtained by decoding the predicted latent variables $\tilde{\bm{y}}(t_n;\boldsymbol{\rho}_p)$. 
The model parameters $\boldsymbol{\theta}_{\enc}$, $\boldsymbol{\theta}_{\dec}$, and $\boldsymbol{\theta}_{\xirom}$ are jointly optimized by minimizing the loss function
\begin{equation}
\begin{split}
      \mathcal{L}_{\text{DL-ROM}}(\boldsymbol{\theta}_{\enc}, \boldsymbol{\theta}_{\dec}, \boldsymbol{\theta}_{\xirom})
      &=\frac{1}{N}\sum_{i=1}^N\Big(\frac{\alpha}{2}\left\|\bm{y}^{(i)} - \dec(\tilde{\bm{y}}^{(i)};\boldsymbol{\theta}_{\dec})\right\|^2_2 + \frac{1-\alpha}{2}\left\|\enc(\bm{y}^{(i)};\boldsymbol{\theta}_{\enc}) - \tilde{\bm{y}}^{(i)}\right\|^2_2 \Big) \\
      &=\frac{1}{N}\sum_{i=1}^N\Big(\frac{\alpha}{2}\left\|\bm{y}^{(i)} - \dec(\xirom(t^{(i)},\boldsymbol{\rho}^{(i)};\boldsymbol{\theta}_{\xirom});\boldsymbol{\theta}_{\dec})\right\|^2_2 + \frac{1-\alpha}{2}\left\|\enc(\bm{y}^{(i)};\boldsymbol{\theta}_{\enc}) - \xirom(t^{(i)},\boldsymbol{\rho}^{(i)};\boldsymbol{\theta}_{\xirom})\right\|^2_2 \Big)\, , 
\end{split}
\label{eq:DL-ROM_training}
\end{equation}
where $0 \le \alpha \le 1$ balances the two terms of the loss function. 
After the optimization of the model parameters (offline phase), in the online phase, the latent dynamics is predicted at any desired $t_n$ and parameter instance $\boldsymbol{\rho}_p$ using $\xirom(t_n,\boldsymbol{\rho}_m;\boldsymbol{\theta}_{\xirom})$, and the corresponding FOM approximations are obtained through the decoder $\dec(\xirom(t^{(i)},\boldsymbol{\rho}^{(i)};\boldsymbol{\theta}_{\xirom});\boldsymbol{\theta}_{\dec})$. 

POD DL-ROM \cite{fresca2022pod} extends DL-ROM by a preliminary POD stage, before the AE neural network, which significantly reduces the dimensionality of the inputs of the AE, and, consequently, the number of trainable parameters and training cost. First, a POD basis of dimension $r$, with $r \ll N_{\bm{y}}$, is computed. Given $U \in \R^{N_{\bm{y}} \times r}$, the POD coefficients are obtained as
\begin{equation}
    \bar{\bm{y}}(t_n;\boldsymbol{\rho}_p)=U^{\top}\bm{y}(t_n;\boldsymbol{\rho}_p).
\end{equation}
The DL-ROM training procedure in \eqref{eq:DL-ROM_training} is then applied to the POD coefficients by replacing the original snapshots with $ \bar{\bm{y}}(t_n;\boldsymbol{\rho}_p)$. In the online phase, latent states are predicted as in DL-ROM and decoded to reconstruct the POD coefficients. Eventually, the corresponding FOM approximations are recovered 
\begin{equation}
    \hat{\bm{y}}(t_n;\boldsymbol{\rho}_m)=U\bar{\bm{y}}(t_n;\boldsymbol{\rho}_m)=U\dec(\xirom(t_n,\boldsymbol{\rho}_m;\boldsymbol{\theta}_{\xirom});\boldsymbol{\theta}_{\dec}).
\end{equation}

\section{Neural network architecture and training}\label{app:neural_network}
\subsection{Autoencoder}
In this work we tested three different architectures for the autoencoder $(\Phi(\cdot;\bm{\theta_{\Phi}}),\Psi(\cdot;\bm{\theta_{\Psi}}))$ to reduce the POD coordinates. They were all trained with the AdamW optimizer, with a learning rate of $10^{-4}$, weight decay of $10^{-4}$ and $\beta_1=0.9,\beta_2=0.98$. Training was done for $1500$ epochs or until no improvement was made on a validation test set for $100$ consecutive epochs. The model with lowest validation loss was selected.
\paragraph{mlp-AE}
This architecture features an encoding stage that takes the POD coordinates as input and has four sequential fully-connected hidden layers of $512$ neurons, each with leaky-ReLU nonlinearity and a final fully-connected layer to bring the dimensionality to $n_{\bm{y}}$.
The decoder stage is symmetrical, also including four fully-connected hidden layers of $512$ neurons with leaky-ReLU nonlinearity and a final fully-connected layer to decode the POD coordinates. The leaky-ReLU is given by

\begin{equation*}
    f_{\text{LeakyReLU}}(x)=\begin{cases}
    \begin{aligned}
        x,\quad &x\geq0\\
        a\cdot x\quad &x < 0.
    \end{aligned}
    \end{cases}
\end{equation*}
With the settings used in the experiment, this architecture has $2.1$M trainable parameters.

\paragraph{conv-AE}
The convolutional encoder uses two convolutional sections, followed by one feed-forward section. The convolutional sections each consist of two convolutional layers, with $32$ output channels, followed by a batch normalization. The feed-forward section has a hidden layer of $512$ neurons followed by a fully-connected layer that maps to the reduced dimension. An ELU nonlinearity is applied after each convolution and fully-connected layer, except the final one. The ELU is defined as

\begin{equation*}
    f_{\text{ELU}}(x)=\begin{cases}
    \begin{aligned}
        x,\quad &x\geq0\\
        a\cdot (e^x-1),\quad &x < 0,
    \end{aligned}
    \end{cases}
\end{equation*}

The decoder stage is similar, though it starts with a fully-connected layer, without hidden layer. Subsequently, there are two convolutional sections, each with two transposed convolutional layers with $32$ output channels, followed by batch normalization. ELU nonlinearities are applied after the fully-connected layer and the transpose convolutions. With the settings used in the experiment, this architecture has $1.1$M trainable parameters.

\paragraph{inv-AE}
The invertible AE uses the approach of \cite{botteghi2026deep}. The encoder network is defined similarly to mlp-AE, as it has five fully-connected hidden layers with $512$ neurons. However, the nonlinearity applied here is GELU, following \cite{botteghi2026deep}. This nonlinearity is given by
\begin{equation*}
    f_{\text{GELU}}(x)=x\Upsilon(x),
\end{equation*}
where $\Upsilon$ is the cumulative distribution of a standard Gaussian. An $n_{\bm{y}}$-dimensional latent representation is achieved by truncating the output to only include the first $n_{\bm{y}}$ values. As is inherent to the invertible AE framework, the decoder is found by inverting the encoder. With the settings used in the experiment, this architecture has $5.3$M trainable parameters.

\subsection{Dynamical model}
\paragraph{Noise-prediction Model for DDP-ROM}
Similarly to \cite{Janner2022PlanningSynthesis, ajay2022conditional} we parameterized the noise-prediction model $\boldsymbol{\epsilon}(\cdot;\boldsymbol{\theta}_{\boldsymbol{\epsilon}})$ with a temporal U-Net. The temporal U-Net is composed of 6 residual blocks and each block is composed of two temporal convolutions, each followed by a Mish nonlinearity. The diffusion step $k$ is encoded with sinusoidal position embedding followed by a fully-connected layer with 256 neurons on Mish nonlinearity. The conditioning input $\bm{c}$ is fed to a separate fully-connected network with a layer of 64 and a layer of 256 hidden neurons and Mish nonlinearity. These embeddings (of dimension 64) are concatenated and provided to the first temporal convolution layer of each block. With the settings used in the experiment, this architecture has $15$M trainable parameters. The noise-prediction model is trained using the Adam optimizer for a total of $10^6$ steps with a learning rate of $2\cdot10^{-4}$ and batch size of $16$. Training went on for $10^6$ training steps, or until the loss evaluated on a validation dataset does not improve for $100$ consecutive evaluations, where evaluations were done every $1000$ training steps. We used $K=200$ diffusion steps, a training horizon $N_t=200$, the probability of removing the conditioning inputs equal to $\text{Bern}(\pi=0.25)$, and low temperature sampling with $\alpha=0.5$ such that the variance in the online phase is $\boldsymbol{\Sigma}_{\text{online}}^k =\alpha \boldsymbol{\Sigma}^k$. 
\paragraph{Dynamical model for POD-DL-ROM}
The neural network used for the dynamical model in the POD-DL-ROM method consists of a three fully-connected layers. After the first two layers, a leaky-ReLU nonlinearity is applied. All hidden layers have $512$ dimensions. In the experiments, the neural network has between $2.67\cdot10^5$ and $2.71\cdot10^5$ trainable parameters. The models were trained with the AdamW optimizer, with a learning rate of $10^{-4}$, weight decay of $10^{-4}$ and $\beta_1=0.9,\beta_2=0.98$. Training was done for $1500$ epochs or until no improvement was made on a validation test set for $100$ consecutive epochs. The model with lowest validation loss was selected.

\end{document}